\input epsf
\magnification1200

\centerline{\bf Right-angled Coxeter groups} 
\centerline{\bf with $n$-dimensional Sierpi\'nski compacta as boundaries.}
\bigskip

\centerline{{Jacek \'Swi\c atkowski}
\footnote{*}{This work was partially supported by the Polish National Science Centre (NCN), grant 2012/06/A/ST1/00259.}}

\medskip
\centerline{Instytut Matematyczny, Uniwersytet Wroc\l awski}
\centerline{pl. Grunwaldzki 2/4, 50-384 Wroc\l aw, Poland}
\centerline{e-mail: {\tt swiatkow@math.uni.wroc.pl}}


\bigskip\noindent
{\bf Abstract.}
For arbitrary positive integer $n$, we describe  
a large class of right-angled Coxeter systems
$(W,S)$ for which the visual boundary $\partial_\infty(W,S)$ is homeomorphic
to the $n$-dimensional Sierpi\'nski compactum.
We also provide a necessary and sufficient condition for a planar simplicial
complex $L$ under which
the right angled Coxeter system $(W,S)$ whose nerve is $L$ has the visual
boundary $\partial_\infty(W,S)$ homeomorphic to the Sierpi\'nski curve.

\bigskip\noindent
{\bf 1. Introduction.}

\medskip
The {\it boundary at infinity} (or the {\it visual boundary}) 
of a Coxeter system $(W,S)$
is a compact metric space canonically associated to the system,
denoted $\partial_\infty(W,S)$. It reflects some aspect of the large-scale
behaviour of the associated Coxeter group $W$, and in case when $W$
is word-hyperbolic, it coincides with the Gromov boundary $\partial W$
of $W$. (We refer the reader to Chapter 12 in [Da]
for a detailed exposition of this concept.)
In this paper, for arbitrary positive integer $n$, we describe,  
a large class of right-angled Coxeter systems
$(W,S)$ for which the boundary $\partial_\infty(W,S)$ is homeomorphic
to the $n$-dimensional Sierpi\'nski compactum.

\medskip\noindent
{\bf 1.1 Definition.}
The {\it $n$-dimensional Sierpi\'nski compactum} is a compact metric space
described, uniquely up to homeomorphism, as any subspace $\Pi$ of the sphere $S^{n+1}$ satisfying the following conditions:
\item{(1)} the complement of $\Pi$ is dense in $S^{n+1}$;
\item{(2)} each connected component $U$ of $S^{n+1}\setminus\Pi$
is an open $(n+1)$-cell in $S^{n+1}$, i.e. the pair $(S^{n+1},U)$ is homeomorphic
to the standard pair $(S^{n+1},\hbox{int}(B^{n+1}))$,
where $B^{n+1}$ denotes here the hemisphere in $S^{n+1}$;
\item{(3)} the family $\cal U$ of all connected components of the 
complement $S^{n+1}\setminus\Pi$ is null, i.e. for any $\epsilon>0$ only finitely many of
these components have diameter greater than $\epsilon$;
\item{(4)} closures $\overline U$ of the components 
$U\in{\cal U}$ are pairwise disjoint.

\medskip
The uniqueness property in the above definition was proved by Cannon [Ca],
see also Theorem 7.2.7 in [DV].
The set of conditions appearing in this definition is called the 
{\it positional characterization}
of the $n$-dimensional Sierpi\'nski compactum.
The 1-dimensional Sierpi\'nski compactum coincides with the space better known
as Sierpi\'nski surve or Sierpi\'n\-ski carpet.

We now explain the special terminology appearing in the statement of the main result of the paper,
Theorem 1.3 below.

\medskip\noindent
{\bf 1.2 Definition} (sphere with holes). An {\it $(n+1)$-sphere with holes} is a simplicial complex
$L$ equipped with a PL embedding into the PL $(n+1)$-sphere $S^{n+1}$
satisfying the following properties:
\item{(0)} $L$ is a proper subspace of $S^{n+1}$;
\item{(1)} $L$ is flag;
\item{(2)} for any connected
component $\Omega$ of $S^{n+1}\setminus L$,
denoting its closure in $S^{n+1}$ by $\overline\Omega$, 
the pair $(S^{n+1},\overline\Omega)$
is PL homeomorphic to the standard pair $(S^{n+1},B^{n+1})$,
where the topological boundary 
$\hbox{bd}(\Omega)=\overline\Omega\setminus\Omega$ corresponds 
to the boundary sphere $\partial B^{n+1}$;
\item{(3)} the boundary $\hbox{bd}(\Omega)$ of any component $\Omega$ of $S^{n+1}\setminus L$ is a full subcomplex of $L$.

\medskip
The term {\it 3-convexity} appearing below in the statement
of Theorem 1.3 is explained in Section 4 
(see Definition 4.1). In the same section we show that for any
positive integer $n$ there are many simplicial complexes $L$
to which the theorem applies.

\medskip
\noindent
{\bf 1.3 Theorem.}
{\it Let $(W,S)$ be a right-angled Coxeter system whose nerve $L$ is an
$(n+1)$-sphere with holes. Suppose also that
 the following two conditions hold:}
\item{(1)} {\it for any two distinct connected components of $S^{n+1}\setminus L$,
their closures in $S^{n+1}$ are either disjoint, or intersect at
a single simplex of $L$;}
\item{(2)} {\it the boundary $\hbox{\rm bd}(\Omega)$ 
of any component $\Omega$ of $S^{n+1}\setminus L$ is a 3-convex 
subcomplex of $L$.}

\noindent
{\it Then the visual boundary $\partial_\infty(W,S)$ is homoemorphic to the $n$-dimensional
Sierpi\'nski compactum.}

\medskip
\centerline{\epsfbox{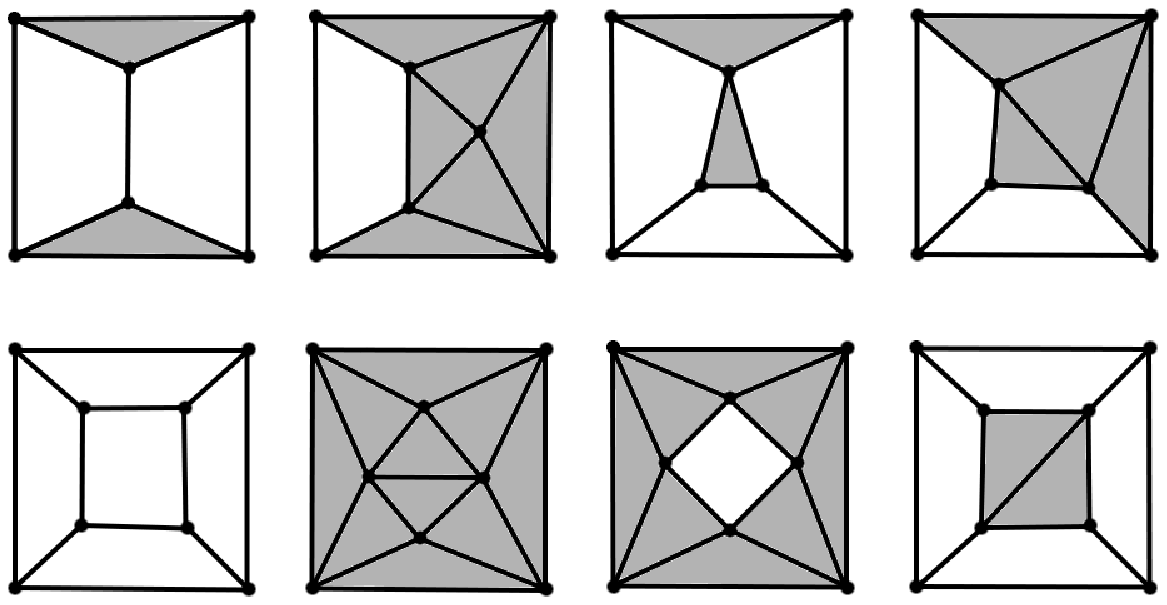}}

\centerline{Figure 1. Examples of nerves of non-hyperbolic right-angled Coxeter groups}
\centerline{with Sierpi\'nski carpet boundaries.}

\bigskip

\medskip
In the case $n=1$, Theorem 1.3 can be used to extend Theorem 1 of [Sw2]
(which characterizes those right-angled word-hyperbolic Coxeter groups with
planar nerves whose Gromov boundary is homeomorphic to the Sierpi\'nski curve)
to non-hyperbolic case.
More precisely, we call
a simplicial complex {\it unseparable} if it is connected, has no separating
simplex, no separating pair of nonadjacent vertices, and no separating
full subcomplex isomorphic to the (simplicial) suspension of a simplex.
Theorem 1.3 yields then the following.

\medskip\noindent
{\bf 1.4 Corollary.}
{\it Let $(W,S)$ be a right angled Coxeter system with planar nerve $L$.
Then the boundary $\partial_\infty(W,S)$
is homeomorphic to the Sierpi\'nski curve if and only if $L$
is unseparable and
distinct from a simplex and from a triangulation of $S^2$.}

\medskip
Examples of some nerves of non-hyperbolic right-angled Coxeter systems with Sierpi\'n\-ski carpet boundaries are presented at Figure 1.
The "only if" part of the corollary is provided by Lemma 1.1 in [Sw2].
Thus,
in order to deduce the corollary from Theorem 1.3,
it is sufficient to note that if a planar flag simplicial complex $L$ 
is unseparable, distinct from a simplex and from a triangulation of $S^2$
then it
is a 2-sphere with holes satisfying conditions (1) and (2) of Theorem 1.3.
We omit a straightforward proof of this observation.

\medskip\noindent
{\bf 1.5 Remark.}
As we show below
(in Lemmas 4.8 and 4.9),
examples of Coxeter systems $(W,S)$ to which Theorem 1.3 applies,
and for which the corresponding group $W$ is word-hyperbolic,
appear only in dimensions $n\in\{ 1,2 \}$.
(Non-hyperbolic examples appear in arbitrary dimension.)
This leads to a question, whether there are any right-angled
word-hyperbolic Coxeter groups with Gromov boundary homeomorphic to
the $n$-dimensional Sierpi\'nski compactum, for $n\ge3$.
This question seems to remain open.

\medskip
The paper is organized as follows. Sections 2 and 6 contain essential
arguments for proving Theorem 1.3, and Sections 3--5 present various
preparatory results and techniques applied in these essential arguments.
More precisely, in Section 2 we present this part of the argument
which uses only the assumption that the nerve $L$
of a right-angled Coxeter system $(W,S)$ is an $(n+1)$-sphere with holes.
Under this restricted assumption, we construct an embedding
of $\partial_\infty(W,S)$ into $S^{n+1}$, and we show that under 
this embedding the subspace $\partial_\infty(W,S)$ satisfies conditions
(1) and (2) of the positional characterization (given in Definition 1.1) 
of the $n$-dimensional Sierpi\'nski compactum.
In Section 6 we use the full strength of the assumptions of Theorem 1.3,
i.e. we assume that the sphere with holes $L$ satisfies also assumptions
(1) and (2) from the statement. We then show that, under the embedding
into $S^{n+1}$ described in Section 2, the subspace $\partial_\infty(W,S)$
satisfies conditions (3) and (4) of Definition 1.1.
This concludes the proof that $\partial_\infty(W,S)$ is homeomorphic to
the $n$-dimensional Sierpi\'nski compactum.

Section 3 contains a preparatory result of independent interest, namely
that the visual boundary of a $\hbox{CAT}(0)$ $(n+1)$-dimensional
PL manifold with nonempty connected convex boundary is homeomorphic to
the $n$-disk $B^n$ (see Proposition 3.1, which is in fact slightly stronger).
This result is essentially used in Section 2, to show that under the described
embedding of $\partial_\infty(W,S)$ in $S^{n+1}$ the components
of the complement are open $(n+1)$-cells in $S^{n+1}$
(which corresponds to condition (2) in Definition 1.1).

In Section 4 we present the concept of 3-convexity, which appears
in assumption (2) of the statement of Theorem 1.3.
We then derive basic properties of this concept, and discuss existence
of nerves $L$ satisfying all the assumptions of Theorem 1.3,
for any $n$.

In Section 5 we consider convex subcomplexes in $\hbox{CAT}(0)$
cubical complexes which satisfy the additional assumption of
{\it local 3-convexity} (which means that the link of the subcomplex at its
any vertex is a 3-convex subcomplex in the corresponding vertex link
of the whole complex). We then deduce some rather strong global geometric
property of such subcomplexes, which resembles the behaviour of
convex subsets in hyperbolic spaces (see Proposition 5.2 for the
statement). This global property plays crucial role in our
arguments in Section 6. The concept of local 3-convexity 
in $\hbox{CAT}(0)$ cubical complexes seems to be
interesting on its own, and potentially useful for other applications.

In Section 6 (and more precisely in Lemma 6.3) we describe one more
technical concept of independent interest, namely a natural family of metrics
(compatible with the topology) on the visual boundary $\partial_\infty X$
of a complete $\hbox{CAT}(0)$ space $X$. Surprisingly, these metrics
seem to be not widely known. We use one of these metrics for showing
that the family of connected components of the complement
$S^{n+1}\setminus\partial_\infty(W,S)$ is null
(which corresponds to condition (3) in Definition 1.1).


\bigskip
\noindent
{\bf 2. Some conclusions for  
arbitrary $(n+1)$-spheres
with holes $L$.}

\medskip
In this section we work under the assumption
that the nerve $L$ of a right-angled Coxeter system $(W,S)$ is an
$(n+1)$-sphere with holes (i.e. we ignore conditions (1) and (2)
from the statement of Theorem 1.3).
We show the following result, which is a first step
in the proof of Theorem 1.3. (Note that assertions (1) and (2) below
correspond to conditions (2) and (1) in Definition 1.1, respectively.)

\medskip\noindent
{\bf 2.1 Proposition.}
{\it Let $(W,S)$ be a right-angled Coxeter system whose nerve $L$ is an
$(n+1)$-sphere with holes. Then}

\item{(1)} {\it there is an embedding of $\partial_\infty(W,S)$ in $S^{n+1}$
such that the connected components of its complement
are open $(n+1)$-cells in $S^{n+1}$;}
\item{(2)} {\it topological dimension of the boundary satisfies
$\dim(\partial_\infty(W,S))=n$, and hence for any embedding of
$\partial_\infty(W,S)$ in $S^{n+1}$ the complement is dense in $S^{n+1}$.}

\medskip
To get assertion (1) above, we first describe and analyze
below, after quick review of the general geometric bockground, a specific embedding
of $\partial_\infty(W,S)$ in $S^{n+1}$. (We will use this embedding
not only here, but also in the
later sections of the paper, in other parts of the proof of Theorem 1.3.) 
We also make essential use of a technical result, Proposition 3.1,
whose statement and proof are presented in the next section.
A short proof of assertion (2) is given at the end of this section.  

\medskip
Throughout this paper we will often work with $CAT(0)$ geodesic metric spaces,
as described e.g. in [BH]. Given  a complete 
CAT(0) space $E$, we denote by $\partial_\infty E$
the visual boundary (at infinity) of $E$,
as defined e.g. in Chapter II.8 in [BH]. We note that if $F\subset E$
is a closed convex subspace of a complete CAT(0) space, then
$F$ is also complete and CAT(0), and the boundary $\partial_\infty F$
is canonically a subset in the boundary $\partial_\infty E$.

Most of the $CAT(0)$ spaces under our interest will be cubical complexes.
Recall that a cubical complex $X$ is CAT(0) (for the standard piecewise euclidean metric $\hbox{d}_X$) iff $X$ is connected, simply connected, and its every
vertex link $\hbox{Lk}(v,X)$ is a flag simplicial complex (cf. Theorem 5.20
on p. 212 in [BH]). A subcomplex $Y$ of a CAT(0) cubical complex $X$
is convex (in the sense that for any two points of $Y$ the geodesic
in $X$ connecting these points is contained in $Y$) iff $Y$ is connected
and for each vertex $v\in Y$ the link $\hbox{Lk}(v,Y)$ is a full
subcomplex in the corresponding link $\hbox{Lk}(v,X)$.

Recall also that for any right-angled Coxeter system $(W,S)$, its nerve $L$
is a flag simplicial complex, and its associated Coxeter-Davis complex $\Sigma$
is a $CAT(0)$ cubical complex whose every vertex link is isomorphic to $L$
(see Proposition 7.3.4 and Theorem 12.2.1(i) in [Da]).
If $W_T$ is a special subgroup of $W$ corresponding to a subset
$T\subset S$, then $(W_T,T)$ is also a right-angled Coxeter system, and its nerve $L_T$
is canonically a full subcomplex of $L$. Moreover, the Coxeter-Davis complex 
$\Sigma_T$ of the system $(W_T,T)$ is canonically a subcomplex in $\Sigma$,
and for any vertex $v$ of $\Sigma_T$ the pair of links
$({\rm Lk}(v,\Sigma),{\rm Lk}(v,\Sigma_T))$ is isomorphic to the pair $(L,L_T)$.
In particular, $\Sigma_T$ is a convex subcomplex of $\Sigma$ and its visual boundary
$\partial_\infty\Sigma_T$ is a subspace in $\partial_\infty\Sigma$.

\medskip
Now, under notation and assumptions of Proposition 2.1, 
we describe a specific embedding
of $\partial_\infty(W,S)$ in $S^{n+1}$.
Given an embedding
of $L$ in $S^{n+1}$ as in the definition of an $(n+1)$-sphere with holes, 
for each connected component $\Omega$
of the complement $S^{n+1}\setminus L$ attach to $L$, along the subcomplex
$\hbox{bd}(\Omega)$ (which is a PL triangulated $n$-sphere), the simplicial cone over this subcomplex.
This yields a flag PL triangulation $L'$ of the sphere $S^{n+1}$ such that 
$L$ appears in $L'$ as a full subcomplex.
Denote by $(W',S')$ the right-angled Coxeter system with nerve $L'$.
Since the boundary of a special subgroup is naturally a subspace in the boundary
of a right-angled Coxeter group, 
we get that $\partial_\infty(W,S)\subset\partial_\infty(W',S')$. Since the latter boundary
is homeomorphic to the sphere $S^{n+1}$ (see e.g. Theorem (3b.2) in [DJ]
or Corollary 1 in [Dr]), it follows that $\partial_\infty(W,S)$ is embedded in $S^{n+1}$.

\medskip
After fixing an embedding $L<L'$ as above,
note that if $L$ is a trangulation of the sphere $S^n$ then both assertions
of Proposition 2.1 follow from Corollary 3.3 (given in the next section).
Thus, in the remaining part of the argument we assume that $L$ is not
a triangulation of $S^n$. In this setting,
we introduce some further terminology and notation which will be used
both in this section and in the remaining part of the paper.
Denote by  $\Sigma=\Sigma(W,S)$ and $\Sigma'=\Sigma(W',S')$
the Coxeter-Davis complexes associated to the respective Coxeter systems.
We then obviously have $S\subset S'$, $W<W'$
(as a special subgroup)
and $\Sigma\subset\Sigma'$ (as a convex subcomplex in a $CAT(0)$ cubical complex).
We also have $\partial_\infty(W,S)=\partial_\infty\Sigma\subset\partial_\infty\Sigma'=S^{n+1}$
as the explicit manifestation of the above mentioned embedding
of $\partial_\infty(W,S)$ in $S^{n+1}$.

For any connected component $\Omega$ of the complement 
$L'\setminus L$,
let $C=\hbox{bd}(\Omega)$ be the subcomplex of $L'$ coinciding with the topological
boundary of $\Omega$ in $L'$. 
It follows from our assumptions that
$C$ is a PL triangulation of $S^n$, and the pair $(L',C)$ is (up to PL homeomorphism)
the standard pair of PL spheres $(S^{n+1},S^n)$.
We call each subcomplex $C$ as above
a {\it peripheral cycle} of $L$ in $L'$.
In the now considered case, when $L$ is not a triangulation of $S^n$,
each peripheral cycle is a proper subcomplex of $L$.

Note that, by condition (3) in Definition 1.2,
any peripheral cycle $C$ is a full subcomplex of $L$.
Consequently, it is also a full subcomplex of $L'$, and thus
it describes a special subgroup of both $W'$ and $W$,
which we denote $W_C$ (and the corresponding induced
Coxeter system by $(W_C,S_C)$). 
Denote by $\Sigma_C=\Sigma(W_C,S_C)$ the Coxeter-Davis complex
of the system $(W_C,S_C)$, and view it canonically as a subcomplex
in $\Sigma'=\Sigma(W',S')$. Since $C$ is full both in $L'$ and $L$, we get that
$\Sigma_C$ is a convex subcomplex of both $\Sigma$ and $\Sigma'$.

Since the pair of nerves $(L',C)$ is the standard pair of 
PL spheres, locally $\Sigma_C$ looks in $\Sigma'$ like an $(n+1)$-dimensional
hyperplane in the euclidean $(n+2)$-space. Since both $\Sigma_C$
and $\Sigma'$ are connected and simply connected (in fact,
they are PL homeomorphic to the euclidean spaces $E^{n+1}$ and $E^{n+2}$,
respectively, see Theorem 10.6.1 in [Da]), it follows that $\Sigma_C$ 
decomposes $\Sigma'$ into the union of two subcomplexes
$H$ such that $H\setminus\Sigma_C$ is 
a connected component of $\Sigma'\setminus\Sigma_C$.
We call these subcomplexes the {\it halfspaces for $\Sigma_C$ in 
$\Sigma'$}.
Note that each such halfspace $H$ is a convex subcomplex of $\Sigma'$.
This follows from the following two observations:
\item{$\bullet$}
a vertex link of $H$ either coincides with the corresponding vertex link of
$\Sigma'$, or these links form a pair isomorphic to the pair $(L',D)$,
where $D$ is one of the (n+1)-disk subcomplexes in $L'$ bounded by $C$;
\item{$\bullet$}
any subcomplex $D$ as above is a full subcomplex of $L'$.

\noindent
In particular, $H$ is a $CAT(0)$ cubical PL $(n+2)$-manifold with convex boundary
$\partial H=\Sigma_C$. 

Next lemma is a direct consequence of
Proposition 3.1 (a technical result whose statement and proof are 
given in the next section).

\medskip\noindent
{\bf 2.2 Lemma.}
{\it If $H_0,H_1$ are the halfspaces for $\Sigma_C$ in $\Sigma'$  then, viewing the boundaries 
$\partial_\infty H_i,\partial_\infty\Sigma_C$ as subspaces in $\partial_\infty\Sigma'\cong S^{n+1}$,
we have:}
\item{(1)} {\it $\partial_\infty H_0\cap\partial_\infty H_1=\partial_\infty\Sigma_C$ and
$\partial_\infty H_0\cup\partial_\infty H_1=\partial_\infty\Sigma'$;}
\item{(2)} {\it $\partial_\infty H_i$ are closed $(n+1)$-cells in $\partial_\infty\Sigma'$ such that 
$\hbox{\rm bd}(\partial_\infty H_i)=\partial_\infty\Sigma_C$, where {\rm bd}
denotes the topological boundary in $\partial_\infty\Sigma'$;}
\item{(3)} {\it $\partial_\infty H_i\setminus\partial_\infty\Sigma_C$ are the connected components of $\partial_\infty\Sigma'\setminus\partial_\infty\Sigma_C$,
and they are open $(n+1)$-cells in $\partial_\infty\Sigma'$;
we also have $\partial_\infty H_i\setminus\partial_\infty\Sigma_C=\partial_\infty\Sigma'
\setminus\partial_\infty H_{i+1}$, where $i+1$ is viewed modulo 2.}

\medskip\noindent
{\bf 2.3 Lemma.}
{\it Let $C$ be a peripheral cycle of $L$ in $L'$.
Then one of the halfspaces obtained by splitting $\Sigma'$ along $\Sigma_C$ contains $\Sigma$, and the other intersects $\Sigma$
only at $\Sigma_C$.}

\medskip\noindent
{\bf Proof:}
We first give an alternative description of the halfspaces for $\Sigma_C$ in $\Sigma'$. Recall that the group $W'$
canonically coincides with the vertex set of the cubical complex
$\Sigma'$, while $S'$ canonically coincides with the vertex set of $L'$.
Let $D_C^+$ be the subcomplex of $L'$ equal to the $(n+1)$-disk bounded by $C$ and containing $L$. Denote also by $D_C^-$ the complementary
$(n+1)$-disk in $L'$ bounded by $C$. Consider the following condition for a reduced expression $g=s_1\dots s_n$ of an element $g\in W'$
(where all $s_i\in S'$):

\smallskip
\item{($*$)} {the first letter $s_i$ not belonging to $S_C$,
if appears at all in the expression, 
corresponds to a vertex in $D_C^+\setminus C$.}

\smallskip 
\noindent
It follows easily from the Tits' characterization of reduced expressions
in Coxeter groups (see Theorem 3.4.3 in [Da]) that property $(*)$
does not depend on the choice of a reduced expression for $g$,
hence it is a condition for $g$.

Denote by ${\cal W}_C^+$ the set of all $g\in W'$ whose reduced
expresions satisfy $(*)$. Denote also by ${\cal W}_C^-$ the analogous
set for which the disk $D_C^+$ in condition $(*)$ is replaced with
$D_C^-$. It is not hard to observe (and we leave it without any further justification) that the two halfspaces for $\Sigma_C$ in $\Sigma'$  are exactly the subcomplexes of $\Sigma'$ spanned
by the sets ${\cal W}_C^+$ and ${\cal W}_C^-$, viewed as subsets
in the vertex set of $\Sigma'$. 
(Here, by the subcomplex spanned by a vertex set we mean the maximal
subcomplex having this set as the vertex set.)
Accordingly, we denote the two halfspaces by $H_C^+$ and $H_C^-$.

Note that obviously we have $W\subset{\cal W}_C^+$ and
$W\cap{\cal W}_C^-=W_C$. Consequently, we get $\Sigma\subset H_C^+$ and $\Sigma\cap H_C^-=\Sigma_C$, which finishes the proof.

\medskip
In consistency with the proof of Lemma 2.3, we denote by $H_C^+$
this halfspace for $\Sigma_C$ in $\Sigma'$ 
which contains $\Sigma$. We also denote by $H_C^-$ the remaining halfspace. Lemmas 2.2 and 2.3 immediately yield the following.

\medskip\noindent
{\bf 2.4 Corollary.}
{\it $\partial_\infty H_C^-\setminus\partial_\infty\Sigma_C=
\partial_\infty\Sigma'\setminus\partial_\infty H_C^+$ is a connected component
of the complement $\partial_\infty\Sigma'\setminus\partial_\infty\Sigma=S^{n+1}
\setminus\partial_\infty(W,S)$.
This component is an open $(n+1)$-cell in $S^{n+1}$,
and its topological boundary coincides with $\partial_\infty\Sigma_C$.}

\medskip
We now pass to describing all other connected components of 
the complement of
$\partial_\infty\Sigma$ in $\partial_\infty\Sigma'$. Note that the group $W$ acts on
$\partial_\infty\Sigma'$ so that it preserves the subspace $\partial_\infty\Sigma$.
In particular, for any $w\in W$ the set 
$\partial_\infty\Sigma'\setminus w\cdot\partial_\infty H_C^+=
\partial_\infty\Sigma'\setminus \partial_\infty(w\cdot H_C^+)$
is a connected component of 
$\partial_\infty\Sigma'\setminus\partial_\infty\Sigma$.
Moreover, one easily observes that 
\item{(1)} components 
$\partial_\infty\Sigma'\setminus \partial_\infty(w_i\cdot H_C^+)$, $i=1,2$,
coincide iff $w_1$ and $w_2$ are in the same coset of $W/W_C$;
\item{(2)} for distinct peripheral cycles $C_1,C_2$ of $L$ in $L'$,
and for any $w_1,w_2\in W$, the components 
$\partial_\infty\Sigma'\setminus \partial_\infty(w_i\cdot H_{C_i}^+)$
are distinct.

\noindent
Due to the above observation (1), we denote by 
$\partial_\infty\Sigma'\setminus \partial_\infty(\bar w\cdot H_C^+)$, where $\bar w\in W/W_C$
is a coset, the space $\partial_\infty\Sigma'\setminus \partial_\infty(w\cdot H_C^+)$
for any $w\in\bar w$.
As a consequence of both above observations (1) and (2), the family
$$
\partial_\infty\Sigma'\setminus \partial_\infty(\bar w\cdot H_C^+):
C \hbox{ is a peripheral cycle of $L$ in $L'$, and }\bar w\in W/W_C
$$
consists of pairwise distinct connected components of 
$\partial_\infty\Sigma'\setminus\partial_\infty\Sigma$. We will show that in fact
there are no other components in this complement. 

\medskip\noindent
{\bf 2.5 Lemma.}
{\it Each connected component of 
$\partial_\infty\Sigma'\setminus\partial_\infty\Sigma$ has a form
$U=\partial_\infty\Sigma'\setminus \partial_\infty(\bar w\cdot H_C^+)$
for some peripheral cycle $C$ of $L$ in $L'$, and for some coset 
$\bar w\in W/W_C$.}

\medskip\noindent
{\bf Proof:}
Start with observing that 
$$
\Sigma=\bigcap_C\bigcap_{\bar w\in W/W_C}\bar w\cdot H_C^+,
$$
where $C$ runs through all peripheral cycles of $L$ in $L'$.
This follows fairly easily from the description of the vertex sets
of the halfspaces $H_C^+$, as given in the proof of Lemma 2.3.
As a consequence, we have
$$
\partial_\infty\Sigma=
\bigcap_C\bigcap_{\bar w\in W/W_C}\partial_\infty(\bar w\cdot H_C^+).
$$
From this we get that
$$
\partial_\infty\Sigma'\setminus\partial_\infty\Sigma=
\partial_\infty\Sigma'\setminus 
\Bigg(
\bigcap_C\bigcap_{\bar w\in W/W_C}\partial_\infty(\bar w\cdot H_C^+)
\Bigg)
=
\bigcup_C\bigcup_{\bar w\in W/W_C}
\bigg(
\partial_\infty\Sigma'\setminus\partial_\infty(\bar w\cdot H_C^+)
\bigg),
$$
which clearly yields the assertion.

\medskip
Denoting by $\hbox{Per}(L,L')$ the set of all peripheral cycles
of $L$ in $L'$, we can summarize the above discussion as follows.

\medskip\noindent
{\bf 2.6 Corollary.}
{\it Let $L$ be an $(n+1)$-sphere with holes
distinct from a triangulation of the sphere $S^n$, 
and let $L<L'\cong S^{n+1}$
be an embedding as described at the beginning of this section. Then}
$${\cal U}= \{ \partial_\infty\Sigma' \setminus\partial_\infty(\bar w\cdot H_C^+):
C\in\hbox{Per}(L,L'), \bar w\in W/W_C\}$$ {\it is the family of all connected components of the complement $\partial_\infty\Sigma'\setminus\partial_\infty\Sigma$.
Each $U\in{\cal U}$
is an open $(n+1)$-cell in $\partial_\infty\Sigma'=S^{n+1}$.
If $U=\partial_\infty\Sigma'\setminus \partial_\infty(\bar w\cdot H_C^+)$
then its closure $\overline U$ and its topological boundary
$\hbox{\rm bd}(U)$ in $\partial_\infty\Sigma'$ can be
expressed as $\overline U=\partial_\infty(\bar w\cdot H_C^-)$
and $\hbox{\rm bd}(U)=\partial_\infty(\bar w\cdot\Sigma_C)$.}

\medskip\noindent
{\bf Proof of Proposition 2.1(1):}
In the case when $L$ is a triangulation of $S^n$,
the assertion follows from Corollary 3.3 (given in the next section). 
In the remaining case, it follows
from Corollary 2.6.

\medskip\noindent
{\bf Proof of Proposition 2.1(2):}
Denote by $\hbox{vcd}(W)$ the virtual cohomological dimension of $W$.
It follows from results of Mike Davis that
$$
\hbox{vcd}(W)=\max\{ k\,:\,\overline{H}^{k-1}(L\setminus\sigma)
\ne0, \hbox{ for some simplex $\sigma$ of $L$, 
or }\overline{H}^{k-1}(L)
\ne0  \}
$$
(see Corollary 8.5.5 in [Da]). Moreover, by the fact that $W$ acts geometrically
on $\Sigma$ and is virtually torsion-free, the pair 
$(\Sigma\cup\partial_\infty(W,S),\partial_\infty(W,S))$ is a ${\cal Z}$-structure
(in the sense of Bestvina described in [Be]) for any
torsion-free finite index subgroup $H<W$.
Since, by Theorem 1.7 of [Be] we then have $\dim\partial_\infty(W,S)=\hbox{cd}(H)-1$
(where $\hbox{cd}(H)$ is the cohomological dimension of $H$),
it follows that $\dim\partial_\infty(W,S)=\hbox{vcd}(W)-1$. Since $L$
is a proper subpolyhedron in $S^{n+1}$, we get 
(from the above formula of Davis) that
$\hbox{vcd}(W)\le n$,
and hence $\dim(\partial_\infty(W,S))\le n$. 

On the other hand, being an $(n+1)$-sphere with holes,
$L$ contains a full subcomplex $C$ which is a PL triangulation
of the sphere $S^n$. Denote by $(W_C,S_C)$ the Coxeter system
for the special subgroup $W_C<W$ induced by $C$.
By Theorem (3b.2)(iii) of [DJ], the boundary $\partial_\infty(W_C,S_C)$
is then homeomorphic to $S^n$, and thus the boundary $\partial_\infty(W,S)$
contains an embedded copy of $S^n$. Hence $\dim(\partial_\infty(W,S))\ge n$,
which completes the proof.


\bigskip\noindent
{\bf 3. Visual boundary of a $CAT(0)$ manifold with convex boundary.}
\nobreak
\medskip
It is known that the visual boundary of a complete $CAT(0)$ PL $(n+1)$-manifold
is homeomorphic to the sphere $S^{n}$ (see [DJ], Theorem (3b.2)(iii)). 
The aim of this section is to extend this result as follows.

\medskip\noindent
{\bf 3.1 Proposition.}
{\it Let $H$ be a complete piecewise euclidean $CAT(0)$ PL-manifold
of dimension $n+1$ with nonempty connected convex boundary $N=\partial H$.
Then the pair of visual boundaries $(\partial_\infty H,\partial_\infty N)$ is homeomorphic
to the pair $(B^n,S^{n-1})$, where $B^n$ is the standard closed 
euclidean $n$-ball
and $S^{n-1}=\partial B^n$.}

\medskip
Before getting to the proof of the proposition, we present its two
immediate corollaries, which are of independent interest.

\medskip\noindent
{\bf 3.2 Corollary.}
{\it Let $M$ be a complete piecewise euclidean PL $(n+1)$-manifold
and $N$ its closed connected 
PL submanifold of dimension $n$. Suppose that $M$ is $CAT(0)$ and $N$ is convex in $M$. Then the pair of visual boundaries
$(\partial_\infty M,\partial_\infty N)$ is homeomorphic to the standard pair of spheres
$(S^n,S^{n-1})$.}

\medskip
Next corollary concerns arbitrary Coxeter groups, not only the
right-angled ones.

\medskip\noindent
{\bf 3.3 Corollary.}
{\it Let $(W',S')$ be a Coxeter system whose nerve $L'$ is a PL triangulation
of the $n$-sphere. Let $L$ be a full subcomplex of $L'$ 
such that the pair $(L',L)$ is PL homeomorphic to the standard pair of spheres
$(S^n,S^{n-1})$. Let $S$ be the vertex set of $L$, and denote by $(W,S)$ the Coxeter system of the special subgroup induced by $S$.
Then the pair of visual boundaries
$(\partial_\infty(W',S'),\partial_\infty(W,S))$ is homeomorphic to the standard pair of spheres
$(S^n,S^{n-1})$.}

\medskip\noindent
{\bf Proof:} Denote by $\Sigma'$ and $\Sigma$ the Coxeter-Davis complexes
of the systems $(W',S')$ and $(W,S)$, respectively.
Viewed naturally as a subcomplex, $\Sigma$ is convex in $\Sigma'$. 
Moreover, 
under assumptions on $L'$ and $L$, $\Sigma'$ is a complete piecewise euclidean
$CAT(0)$ PL-manifold of dimension $n+1$, and $\Sigma$ is its convex
closed PL submanifold of dimension $n$. Corollary 3.2 yields then
the assertion.

\bigskip
The rest of this section is devoted to the proof of Proposition 3.1.
In this proof we follow the lines of the proof of Theorem (3b.2)(iii) in [DJ]
(see also the appendix in [Sw1]).
In particular, we
use the technique of cellular sets and cell-like maps, whose basic elements we
now recall. 

In our setting below, a manifold is a (not necessarily compact) topological manifold
with boundary (possibly empty). A compact subset $E$ of a manifold $M$ is 
{\it cellular} if the quotient space $M/E$ (obtained by collapsing $E$ to a point) is homeomorphic to $M$, and the quotient map $q:M\to M/E$
can be approximated by homeomorphisms. The latter means that
for any fixed metric $d$ in $M/E$ (compatible with the topology)
and any continuous positive function $\delta:M\to R_+$
there is a homeomorphism $h:M\to M/E$ such that for all $x\in M$
we have $d(q(x),h(x))<\delta(x)$. 

Cellular subsets of a manifold can be also characterised as admitting
systems of neighbourhoods of a special form. To make a precise statement, recall that an {\it $n$-cell} in an $n$-dimensional manifold $M$
is any subset $B\subset M$ which in some neighbourhood in $M$ looks like:
\item{(1)} $B^n$ in $R^n$, if $B\subset\hbox{int}(M)$;
\item{(2)} $B^n_+$ in $R^n_+$, if $B\cap\partial M\ne\emptyset$,
where $R^n_+=\{ (x_1,\dots,x_n):x_n\ge0 \}$ and 
$B^n_+=\{ x\in R^n_+:|x|\le1 \}$.

\noindent
An {\it open $n$-cell}\hskip6pt in $M$ is a subset of an $n$-cell corresponding
to the open ball $(B^n)^\circ=B^n\setminus\partial B^n$ or the open half-ball
$(B^n_+)^\circ=\{ x\in R^n_+:|x|<1 \}$, respectively.

\medskip\noindent
{\bf 3.4 Lemma} (see Theorem 2.2 in [He]).
{\it A compact subset $E$ in an $n$-dimensional manifold $M$ is cellular iff
it admits arbitrarily close neighbourhoods in $M$ which are open $n$-cells.
More precisely, fixing any metric $d$ in $M$ (compatible with the topology),
for any $\varepsilon>0$ there is an open $n$-cell $U$ in $M$ such that
$E\subset U\subset N_\varepsilon(E)$, where $N_\varepsilon(E)=\{ x\in M :d(x,E)<\varepsilon \}$.}

\medskip
A nonempty compact metric space of finite topological dimension
is {\it cell-like} if it can be
embedded as a cellular subset in some manifold
(for justification of this characterization, see the remark in the middle
of page 114 in [Ed]). Obviously,
each cellular subset of a manifold is cell-like, but the converse
is not true (i.e. there are non-cellular embeddings of cell-like spaces
in manifolds). A map between metric spaces is {\it cell-like}
if it is a proper surjection and each point preimage is cell-like.

\medskip
We now recall few preparatory results which will be needed in our proof
of Proposition 3.1.
The next result is due to Siebenmann [Si] for $n\ge5$,
Quinn [Qu] for $n=4$, Armentrout [Ar] for $n=3$ and R. L. Moore
[Mo] for $n\le2$. The case of manifolds with nonempty boundary
is carefully addressed in [Si], and the cases involving dimension 3 hold
because the Poincare conjecture is true.

\medskip\noindent
{\bf 3.5 Approximation Theorem.}
{\it Each cell-like map between manifolds with boundary 
is a near-homeomor\-phism,
i.e. it can be approximated by homeomorphisms. More precisely,
if $f:M\to N$ is cell-like and $d$ is any metric in $N$ (compatible
with the topology) then for any continuous positive function
$\delta:M\to R_+$ there is a homeomorphism $h:M\to N$
such that for all $x\in M$ we have $d(f(x),h(x))\le\delta(x)$.}

\medskip
The next result seems to be well known 
(it appears as Lemma 2 in the appendix to [Sw1], where its proof is included). 

\medskip\noindent
{\bf 3.6 Lemma.}
{\it If $f:A\to B$ is a cell-like surjection between finite dimensional
metric compacta, then $A$ is cell-like if and only if $B$ is cell-like.}

\medskip
Next two preparatory results give some useful criteria for recognizing cellular
subsets in spheres and disks. First of them presents a well known
characterization of cellular subsets
in spheres $S^n$ as those which are {\it point-like} (i.e. such that the complement
is homeomorphic to $R^n$). The second one is a less well known analogue
for those subsets in disks which intersect the boundary.

\medskip\noindent
{\bf 3.7 Lemma} (see Theorem on page 114 in [Ed]).
{\it A closed subset $E$ of the sphere $S^n$ is cellular iff the complement 
$S^n\setminus E$ is homeomorhic to $R^n$.}

\medskip\noindent
{\bf 3.8 Lemma.}
{\it Let $E$ be a closed subset of the $n$-disk $D^n$ such that 
$E\cap\partial D^n\ne\emptyset$. Suppose that the complement 
$D^n\setminus E$ is homeomorphic to $R^n_+$. Then $E$ is cellular in $D^n$.}

\medskip\noindent
{\bf Proof:}
We need to indicate arbitrarily close open neighbourhoods of $E$ in $D^n$
which are open $n$-cells (of the form $(B^n_+)^\circ$).
Let $h:R^n_+\to D^n\setminus E$ be a hypothesized homeomorphism, and
let $B_r$ denote the euclidean closed ball of radius $r$ in $R^n_+$
centered at $0$. Obviously,   as $r\to\infty$, the set $D^n\setminus h(B_r)$
becomes as close  to $E$ as necessary.
We claim that
$D^n\setminus h(B_r)$ is also an open $n$-cell in $D^n$.
This follows fairly directly from the generalized Schoenflies theorem.

\medskip
Having finished the preparations,
we now state two auxilliary results which will lead us to the proof of
Proposition 3.1. The term $\pi$-{\it convexity} appearing in Lemma 2
means that any geodesic of length less than $\pi$ having both endpoints in
the subspace, is entirely contained in this subspace.

\medskip\noindent
{\bf Lemma 1.}
{\it Let $n$ be a positive integer, and 
let $H$ be a complete piecewise euclidean $CAT(0)$ PL-manifold
of dimension $n$, with nonempty connected convex boundary $N=\partial H$. 
Then for any $x\in N$ and any $R>0$:} 
\itemitem{$(B(n))$}
{\it the sphere $S_R(x,H)$ is an $(n-1)$-manifold and $S_R(x,N)$ is its
boundary;}
\itemitem{$(C(n))$}
{\it for any $R'>R$ the geodesic projection towards $x$,
$\rho_{R',R}:S_{R'}(x,H)\to S_R(x,H)$, is a cell-like map;}
\itemitem{$(D(n))$}
{\it the pair  $(S_R(x,H),S_R(x,N))$ is homeomorphic to the pair $(B^{n-1},S^{n-2})$.}

\medskip\noindent
{\bf Lemma 2.}
{\it Let $n$ be a positive integer, and let
$P$ be a complete piecewise spherical $CAT(1)$ PL-manifold
of dimension $n$, with nonempty 
$\pi$-convex boundary $Q=\partial P$.
Then for any $v\in Q$:} 
\itemitem{$(A_S(n))$}
{\it for any $R\in(0,\pi]$ the open metric ball $B_R(v,P)$ is homeomorphic to
the open $n$-dimensional half disk, and $B_R(v,Q)$ corresponds to its
boundary;}

\itemitem{$(C_S(n))$}
{\it for any $R<R'<\pi$ the geodesic projection towards $v$,
$\rho_{R',R}:S_{R'}(v,P)\to S_R(v,P)$,  is a cell-like map.}

\medskip
We will prove all the assertions of Lemmas 1 and 2 by induction with
respect to the dimension $n$. For $n=1$ all of these statements are
obviously true. The main inductive step will consist of the following implications:
\item{$\bullet$} $A_S(n-1)\Rightarrow [B(n)\wedge C_S(n)
\wedge C(n)]$;
\item{$\bullet$} $[B(n)\wedge C(n)]\Rightarrow D(n)$;
\item{$\bullet$} $C_S(n)\Rightarrow A_S(n)$.

\noindent
We will also show (prior to proving the lemmas) the implication
$(\forall n \,\,C(n)\wedge D(n))\Rightarrow$ Proposition 3.1.

\medskip\noindent
{\bf Proof of Proposition 3.1} (assuming $\forall n \,\,C(n)\wedge D(n)$):

Let $H$ and $N$ satisfy the assumptions of the proposition
(in particular, $\dim H=n+1$), and let $x\in N\subset H$ be any point. 
Recall that we have 
$$
\partial_\infty H=\lim_{\longleftarrow}[\{ S_R(x,H):R>0 \}, \{ \rho_{R',R}:R'>R \}]
$$
and
$$
\partial_\infty N=\lim_{\longleftarrow}[\{ S_R(x,N):R>0 \}, 
\{ \rho_{R',R}|_{S_{R'}(x,N)}:R'>R \}],
$$
where $\rho_{R',R}$ are the geodesic projections. 
Since by $C(n+1)$ these bonding maps $\rho_{R',R}$ are all cell-like,
it follows from Approximation Theorem 3.5 that they are near-homeomor\-phisms. By a result of M. Brown (Theorem 4 in [Br]),
if we fix an increasing sequence of radii $R_i\to\infty$,
then we can also choose homeomorphisms $h_i:S_{R_i}(x,H)\to S_{R_{i-1}}(x,H)$ (appropriately close to the maps $\rho_{R_{i},R_{i-1}}$) such that
the visual boundary $\partial_\infty H$ is homeomorphic to the inverse limit
$$
\lim_{\longleftarrow}[\{ S_{R_i}(x,H) \},\{ h_i \}].
$$
By $D(n+1)$, this inverse limit is homeomorphic to $B^{n}$.
Moreover, according to Brown's description of
an explicit homeomorphism 
$F:\partial_\infty H\to\lim_{\longleftarrow}[\{ S_{R_i}(x,H) \},\{ h_i \}]$
as above (provided in the proof of Theorem 1 in [Br]), any point
$z\in\partial_\infty H$ which is represented by a string $(z_i)$ 
(of the inverse system $[\{ S_{R_i}(x,H) \}, \{ \rho_{R_i,R_{i-1}} \}]$)
is mapped by $F$ to a point $y$ represented by a string $(y_i)$
given by $y_i=\lim_{j>i}h_{i}\circ\dots\circ h_{j-1}(z_j)$.
In particular, if $z\in\partial_\infty N$, we have that all  $z_i\in S_{R_i}(x,N)$,
and consequently all $y_i\in S_{R_i}(x,N)$ too,
because homeomorphisms preserve boundary. 
It follows that under $F$ the subset $\partial_\infty N\subset\partial_\infty H$
is mapped to the subset 
$$
\lim_{\longleftarrow}[\{ S_{R_i}(x,N) \},\{ h_i|_{S_{R_i}(x,N)} \}]
\cong S^{n-1}\subset B^n\cong \lim_{\longleftarrow}[\{ S_{R_i}(x,H) \},\{ h_i \}].
$$
Since we know that $\partial_\infty N$ is homeomorphic to $S^{n-1}$
(by the result from [DJ] mentioned at the beginning of this section),
it follows that $F$ maps $\partial_\infty N$ onto the above subset,
and this completes the proof.

\medskip
We now turn to the inductive proof of Lemmas 1 and 2,
following the scheme outlined right after the statement of Lemma 2.

\medskip\noindent
{\bf Proof of the implication $A_S(n-1)\Rightarrow  B(n)$:}

We refer to the argument as in the proof of the implication 
$(L_{n-1})\Rightarrow (T_n)$ in [DJ], in a long paragraph beginning at the end
of page 371 and ending in the middle of the next page.
This argument, without any change, shows that $S_R(x,H)\setminus S_R(x,N)$
is an $(n-1)$-manifold. Thus, to get $B(n)$, it is sufficient to show that
any $y\in S_R(x,N)$ has an open neighbourhood in $S_R(x,H)$ homeomorphic
to the $(n-1)$-dimensional open half-disk $(B^{n-1}_+)^\circ$.
For the latter, it is obviously sufficient to show that any $y$ as above
has a closed neighbourhood $\Delta$ in the closed ball $\overline B_R(x,H)$
such that the pair $(\Delta,\Delta\cap \overline B_R(x,N))$  is homeomorphic
to the pair $(B^n_+,B^{n-1}_0)$, with $y$ corresponding to a point
on the boundary $\partial B^{n-1}_0$.

So let $y\in S_R(x,N)$ be any point. 
Denote by $v\in\hbox{Lk}(y,H)$ the point induced by the geodesic $[y,x]$.
The argument from [DJ] mentioned above,
after applying to our setting, shows that for small $\varepsilon>0$ the ball
of radius $\varepsilon$ in $\overline B_R(x,H)$ centered at $y$
is homeomorphic to the quotient 
$$
\overline B_{\pi/2}(v,\hbox{Lk}(y,H))\times[0,\varepsilon] \,\,/ \,\,
\overline B_{\pi/2}(v,\hbox{Lk}(y,H))\times\{0\}\cup
S_{\pi/2}(v,\hbox{Lk}(y,H))\times[0,\varepsilon],
$$  
where the intersection of this $\varepsilon$-ball with $B_R(x,N)$
corresponds to the subset 
$$
\overline B_{\pi/2}(v,\hbox{Lk}(y,N))\times[0,\varepsilon] \,\,/ \,\,
\overline B_{\pi/2}(v,\hbox{Lk}(y,N))\times\{0\}\cup
S_{\pi/2}(v,\hbox{Lk}(y,N))\times[0,\varepsilon]
$$
and where $y$ corresponds to the point obtained from the collapsed set.
Now, it follows from $A_S(n-1)$  
that the pair of open metric balls
$(B_{\pi/2}(v,\hbox{Lk}(y,H)), B_{\pi/2}(v,\hbox{Lk}(y,N)))$ is homeomorphic
to the standard pair of open balls $((B^{n-1}_+)^\circ, (B^{n-2}_0)^\circ)$.
This easily implies the last assertion of the previous paragraph,
and hence completes the proof.

\medskip\noindent
{\bf Proof of the implication $A_S(n-1)\Rightarrow[C_S(n)\wedge C(n)]$:}

We present a proof of the implication $A_S(n-1)\Rightarrow C(n)$,
and it will become clear that the same argument yields also
the implication $A_S(n-1)\Rightarrow C_S(n)$.

Let $p\in S_R(x,H)$, and denote shortly by $\rho $ the geodesic projection 
$\rho_{R',R}:S_{R'}(x,H)\to S_R(x,H)$.
To prove that $\rho$ is cell-like, we need to show that 
the preimage $\rho^{-1}(p)$ is a cell-like set.
We will show this for the case when $p\in S_R(x,N)$, by an argument following
that in the proof of Lemma (3b.1)
and Theorem (3b.2) in [DJ] (see also the appendix to [Sw1]).
The argument in the remaining case $p\in S_R(x,H)\setminus S_R(x,N)$
is analogous, follows more directly from that in [DJ],
builds upon Lemma 3.7 rather than Lemma 3.8, and we omit its details.

So suppose $p\in S_R(x,N)$. Our proof consists of two steps which
correspond to the following two claims.

\smallskip\noindent
{\bf Claim 1.}
{\it Let $y\in N=\partial H$, and 
let $v$ be any point from the sublink 
$\hbox{\rm Lk}(y,N)\subset \hbox{\rm Lk}(y,H)$
(where the latter pair of links is homeomorphic to the standard pair
$(B^{n-1},\partial B^{n-1})$).
Then the complement in $L:=\hbox{\rm Lk}(y,H)$ of the open ball $B_\pi(v,L)$
is a cellular subset of $L$.}

\smallskip
To prove Claim 1,
denote the complement $L\setminus B_\pi(v,L)$ by $\hbox{Sh}_L(v)$
and call it the {\it shadow} of $v$ in $L$. 
Since $N$ is a complete $CAT(0)$ manifold, it satisfies the geodesic extension
property, and since it is convex in $H$, it follows that the shadow 
$\hbox{Sh}_L(v)$ has nonempty intersection with the sublink 
$\hbox{\rm Lk}(y,N)$. Thus, $\hbox{Sh}_L(v)$ is a closed subset
in the $(n-1)$-disk $L$, and it intersects the boundary $\partial L$.
On the other hand, by $A_S(n-1)$, the complement 
$L\setminus\hbox{Sh}_L(v)=B_\pi(v,L)$ is homeomorphic to the open
$(n-1)$-dimensional half-disk. It follows then from Lemma 3.8 that
$\hbox{Sh}(v,L)$ is cellular in $L$, hence the claim.

\smallskip\noindent
{\bf Claim 2.}
{\it If $p\in S_R(x,N)$ then the preimage $\rho^{-1}(p)$ is cell-like.}

\smallskip
To prove Claim 2 we use an idea communicated to us by Mike Davis.
Suppose {\it a~contrario} that the preimage $\rho^{-1}(p)$ is not cell-like.
Consider all geodesic rays in $H$ started at $x$ and passing through $p$,
and let $[x,\tilde p]$ be the largest geodesic segment contained in all of those
geodesic rays (it may happen that $\tilde p=p$). Consider the link 
$\hbox{\rm Lk}(\tilde p,H)$
and the point $\tilde v\in \hbox{\rm Lk}(\tilde p,H)$ 
corresponding to the geodesic segment $[x,\tilde p]$.
Then the shadow of $\tilde v$ in $\hbox{\rm Lk}(\tilde p,H)$ is nontrivial.
Put $\tilde R=d_H(x,\tilde p)$. Since $\tilde p$ has a metric cone neighbourhood in $H$ (see e.g. Theorem 7.16 on page 104 in [BH]), 
there is $\tilde R'>\tilde R$ such that
the preimage $\rho_{\tilde R',R}^{-1}(p)$ is homeomorphic to the shadow $\hbox{Sh}_{{\rm Lk}(\tilde p,H)}(\tilde v)$. 
In particular, it follows from Claim 1 that
$\rho_{\tilde R',R}^{-1}(p)$ is homeomorphic to a cellular subset of 
$\hbox{\rm Lk}(\tilde p,H)$,
and hence it is cell-like.
Since we have $\rho^{-1}(p)=
\rho_{R',\tilde R'}^{-1}[\rho_{\tilde R',R}^{-1}(p)]$,
it follows from Lemma 3.6 that there is $p_1\in \rho_{\tilde R',R}^{-1}(p)$ 
such that $\rho_{R',\tilde R'}^{-1}(p_1)$ is not cell-like.

Iterating the argument as above, we get an infinite sequence of points $p_n$
lying on a common finite geodesic segment started at $x$ with the following
property: denoting by $v_n\in \hbox{\rm Lk}(p_n,H)$ the point corresponding to the segment $[x,p_n]$, we have that the shadow $\hbox{Sh}_{{\rm Lk}(p_n,H)}(v_n)$
is nontrivial. On the other hand, it is not hard to realize that for any geodesic
segment $[x,y]$ the set of points $p\in[x,y]$ with such nontrivial shadows
is discrete, and hence finite (cf. Fact 4 in the appendix to [Sw1]). 
This contradiction concludes the proof of Claim 2,
and hence also the proof
of the implication $A_S(n-1)\Rightarrow C(n)$.

\medskip\noindent
{\bf Proof of the implication $[B(n)\wedge C(n)]\Rightarrow D(n)$:}

Note that, if $r>0$ is sufficiently small, the pair $(S_r(x,H),S_r(x,N))$ is homeomorphic to the pair $(B^{n-1},S^{n-2})$ by the fact that
$x$ has a metric cone neighbourhood in $H$. Since by $C(n)$, 
for any $R>r$
the map $\rho_{R,r}:S_R(x,H)\to S_r(x,H)$ is cell-like,
the assertion of $D(n)$ follows from
Approximation Theorem 3.5 combined with $B(n)$.

\medskip\noindent
{\bf Proof of the implication $C_S(n)\Rightarrow A_S(n)$:}

Note that for small $R$ the assertion of $A_S(n)$ follows from the existence
of a metric cone neighbourhood of $v$ in $P$. Fix some $r>0$ such that
$B_r(v,P)$ is a cone neighbourhood. It remains to prove the assertion for
$R\in(r,\pi]$. Given any $R$ as above, put $c:={r\over R}$, and consider the map 
$f:B_R(v,P)\to B_r(v,P)$ defined as follows. For any $p\in B_R(v,P)$ put
$r_p:=d_P(v,p)$, and then put 
$f(p):=\rho_{r_p,cr_p}(p)$.
Obviously, $f$ is a continuous surjection between open balls
(which are $n$-manifolds), and we claim that $f$ is cell-like. 
Indeed, obviously we have $f^{-1}(v)=\{ v \}$, and
for any $q\in B_r(v,P)$, $q\ne v$, putting $r_q=d_P(v,q)$, the preimage
$f^{-1}(q)$ is easily seen to be homeomorphic to the set 
$\rho_{r_q/c,r_q}^{-1}(q)$. Since $r_q/c<R\le \pi$, it follows from $C_S(n)$
that the latter set is cell-like, and hence $f$ is cell-like as well.
 Now, applying Approximation Theorem 3.5 we get a homeomorphism
 $h:B_R(v,P)\to B_r(v,P)$, which shows that $B_R(v,P)$ is homeomorphic to
 the open $n$-dimensional half-disk.
 Since obviously we have $\partial B_R(v,P)=B_R(v,Q)$, the second part of the
 assertion of $A_S(n)$ follows as well, which completes the proof.


\bigskip
\noindent
{\bf 4. 3-convexity in simplicial complexes.}

\medskip
In this section we present the concept of 3-convexity for subcomplexes
in simplicial complexes, and discuss existence of various classes of
examples of simplicial complexes $L$ satisfying all the assumptions
of Theorem 1.3. 

\medskip\noindent
{\it 3-convexity.}

\smallskip
Given a simplicial complex $K$, we consider 
the combinatorial distance function 
$\hbox{dist}_K$ between its vertices given by minimizing
the length (i.e. the number of edges) of polygonal paths in the 1-skeleton 
$K^{(1)}$ of $K$ connecting these vertices.

\medskip\noindent
{\bf 4.1 Definition} (3-convexity).
A subcomplex $J$ in a simplicial complex $K$ is {\it 3-convex} if it is a full
subcomplex and if for any two vertices of $J$ staying at distance 2 in $K$
any path of length 2 in $K^{(1)}$ connecting these vertices is entirely
contained in $J$.

\medskip
Next lemma presents a useful property of 3-convex subcomplexes.
Given a simplex $\sigma$ in a simplicial complex $K$, we denote by 
$\hbox{st}_K(\sigma)$, and call the {\it star of $\sigma$ in $K$},
the subcomplex of $K$ being the union of all simplices of $K$ containing 
$\sigma$. 

\medskip\noindent
{\bf 4.2 Lemma.}
{\it A full subcomplex $J$ of a flag simplicial complex $K$ is 3-convex
if and only if for any vertex $v\in J$ the union $J\cup\hbox{st}_K(v)$
is a full subcomplex of $K$. Moreover, if $J$ is 3-convex in $K$ then
for any simplex $\sigma\in J$ the union $J\cup\hbox{st}_K(\sigma)$
is a full subcomplex of $K$.}

\medskip\noindent
{\bf Proof:} Let $\sigma\in J$ and let $v_1,\dots,v_m$ be any vertices of $J\cup\hbox{st}_K(\sigma)$ pairwise connected with edges in $K$.
By flagness of $K$, there is a simplex spanned on these vertices,
and we denote it $\tau$. We need to show that 
$\tau\subset J\cup\hbox{st}_K(\sigma)$.

If $v_1,\dots,v_m\in J$, we get $\tau\subset J$ by fullness of $J$,
and this clearly yields the assertion.
We will show that otherwise all vertices $v_1,\dots,v_m$ are contained
in $\hbox{st}_K(\sigma)$. Let $v\in\{ v_1,\dots,v_m \}$ be a vertex
not contained in $J$, and suppose that $u\in\{ v_1,\dots,v_m \}$
is a vertex not contained in $\hbox{st}_K(\sigma)$. Then $u\in J$,
and there is $w\in\sigma$ such that $\hbox{dist}_K(u,w)=2$.
By 3-convexity of $J$, the geodesic path $(u,v,w)$ is then contained in $J$,
which contradicts the choice of $v$. Thus, 
$\{ v_1,\dots,v_m \}\subset\hbox{st}_K(\sigma)$ in this case,
and the assertion follows  by the fact that star of any simplex
in a flag simplicial complex is its full subcomplex.

To complete the proof, it remains to show that if for each vertex $v\in J$
the union $J\cup\hbox{st}_K(v)$ is full in $K$ then $J$ is 3-convex in $K$.
Let $v_1,v_2$ be any vertices of $J$ lying at distance 2 in $K$,
and let $(v_1,w,v_2)$ be any geodesic between them in $K^{(1)}$.
If $w\notin J$ then we have 
$[w,v_2]\not\subset J\cup\hbox{st}_K(v_1)$,
contradicting fullnes of the latter union.

\medskip\noindent
{\it Constructing 3-convex subcomplexes.}

\smallskip
In this subsection, given any polyhedron $P$ and its subpolyhedron $Q$, 
we present a construction
of a flag triangulation of $P$ for which $Q$
forms a 3-convex subcomplex. This construction allows to justify the existence
of many simplicial complexes $L$ satisfying the assumptions
of Theorem 1.3, in arbitrary dimension $n$ (see Example 4.6 below).

\medskip\noindent
{\bf 4.3 Definition} (relative barycentric subdivision).
Let $K$ be a simplicial complex and $J$ its subcomplex.
The {\it barycentric subdivision of $K$ relative to $J$},
denoted $K'_J$, is the subdivision of $K$ obtained as follows.
First, subdivide barycentrically each 1-simplex of $K$ not contained in $J$,
and then, recursively with respect to the dimension, subdivide each
higher dimensional simplex not contained in $J$ by coning its already
subdivided boundary. Alternatively, take as the vertex set of $K'_J$
all vertices of $J$ and the barycenters $b(\sigma)$ of all simplices $\sigma$ of $K$
not contained in $J$; simplices of $K'_J$ are spanned precisely by the sets
of form $\{ v_1,\dots,v_q,b(\sigma_1),\dots,b(\sigma_p) \}$, where $q\ge0$, 
$p\ge0$, the vertices $v_1,\dots,v_q$ span a simplex $\tau$ of $J$ (empty when $q=0$),
$\sigma_1,\dots,\sigma_p$ are some simplices of $K$ not contained in $J$,
and $\tau\subset\sigma_1\subset\dots\subset\sigma_p$.

\medskip
Note that the subcomplex $J\subset K$, with its original simplicial structure induced
from $K$, is also a subcomplex of $K'_J$. Next few lemmas collect various further 
properties of relative barycentric subdivision.

\medskip\noindent
{\bf 4.4 Lemma.} {\sl If a subcomplex $J$ of a simplicial complex $K$ is flag
then the subdivision $K'_J$ is also flag.
In particular, if $J$ is a full subcomplex of a flag complex $K$
then the subdivision $K'_J$ is also flag.}

\medskip\noindent
{\bf Proof:} Let $A$ be a set of vertices of $K'_J$ which are pairwise connected
with edges. Suppose $A=\{ u_1,\dots,u_k,b(\rho_1),\dots,b(\rho_m) \}$, where
$u_1,\dots,u_k$ are some vertices of $J$ and $\rho_1,\dots,\rho_m$ are some simplices
of $K$ not contained in $J$. By the assumption on $A$, the vertices $u_1,\dots,u_k$
are pairwise connected with edges of $J$, and since $J$ is flag, they span a simplex
$\tau$ in $J$. Again, by the assumption on $A$, for any two simplices $\rho_i,\rho_j$
one of them is a face of the other. Consequently, we may assume that 
$\rho_1\subset\dots\subset\rho_m$. Moreover, the existence of edges in $K'_J$
between $b(\rho_1)$ and any of the $u_i$ implies that $\tau\subset\rho_1$.
It follows then directly from the description of simplices in $K'_J$
(in Definition 4.3) that the set $A$ spans a simplex of $K'_J$. 
Hence flagness of $K'_J$.

\medskip\noindent
{\bf 4.5 Lemma.} {\sl If a subcomplex $J$ of a simplicial complex $K$ is flag and full in $K$
then $J$ is a 3-convex subcomplex of $K'_J$.
In particular, if $J$ is a full subcomplex of a flag complex $K$
then $J$ is 3-convex in $K'_J$.}

\medskip\noindent
{\bf Proof:}
The fact that $J$ is a full subcomplex of $K'_J$ is obvious.
To show that $J$ is also 3-convex, consider any vertices $v,w$ of $J$
at combinatorial distance 2 in $K'_J$.
Suppose on the contrary that there is a combinatorial geodesic 
$(v,u,w)$ in the 1-skeleton of $K'_J$ such that $u\notin J$.
By definition of $K'_J$, the 1-simplex $(v,u)$ must have the form
$(v,b_\sigma)$ for some simplex $\sigma$ of $K$ containing $v$,
and $(u,w)$ must have the form
$(b_\tau,w)$ for some simplex $\tau$ of $K$ containing $w$.
It follows that $\sigma=\tau$, and we get that $u$ and $w$
belong to the same simplex of $K$. Since $J$ is full in $K$,
the edge $(v,w)$ is an edge of $J$, and hence also an edge of $K'_J$,
contradicting the fact that $\hbox{dist}_{K'_J}(v,w)=2$.

\medskip\noindent
{\bf 4.6 Example.}
Using Lemmas 4.4 and 4.5, one can construct many simplicial complexes $L$
satisfying the assumptions of Theorem 1.3, in arbitrary dimension $n$.
Indeed, consider any PL triangulation $N$ of the $(n+1)$-sphere,
and any collection $\sigma_1,\dots,\sigma_k$ of pairwise nonintersecting
$(n+1)$-simplices of $N$, with $k\ge1$.
Let $M$ be the simplicial complex obtained from $N$ be deleting
interiors of all simplices $\sigma_1,\dots,\sigma_k$.
Obviously, $M$ is a PL triangulation of an $(n+1)$-manifold
with boundary ($k$ times punctured $(n+1)$-sphere).
Let $M'$ be the first barycentric subdivision of $M$, and let $J$
be its subcomplex corresponding to the boundary.
Then obviously $M'$ is flag, and $J$ is a full subcomplex of $M'$.
Putting $L=(M')'_J$, we obtain (in view of Lemma 4.4) 
a flag simplicial complex for which $J$ is (by Lemma 4.5)
a 3-convex subcomplex. Consequently, the connected components of $J$
(which correspond to subcomplexes $\hbox{bd}(\Omega)$ 
of Definition 1.2) are also 3-convex in $L$.
This shows that $L$ satisfies all the 
assumptions of Theorem 1.3.

Similarly, given any polyhedron $P$ and its any subpolyhedron $Q$,
and starting from any triangulation $M$ of $P$ for which $Q$ induces
a subcomplex, we can apply the same procedure as above to get
a flag triangulation $L$ of $P$ for which $Q$ induces a 3-convex
subcomplex.

\medskip\noindent
{\it Flag-no-square spheres with holes having all peripheral cycles 3-convex.}

\smallskip
In this subsection we discuss existence of examples of $(n+1)$-spheres
with holes, in various dimensions $n$, satisfying all the assumptions
of Theorem 1.3, and such that the associated Coxeter group $W$
is word-hyperbolic.

Recall that,
given an $(n+1)$-sphere with holes $L$ (viewed as PL embedded
in the sphere $S^{n+1}$), a {\it peripheral cycle} in $L$ is any
subcomplex $C$ which is induced by the topological boundary
$\hbox{bd}(\Omega)$ of a connected component $\Omega$
of the complement $S^{n+1}\setminus L$.

A simplicial complex $K$ is said to satisfy the {\it flag-no-square} condition
if it is flag and the following no-empty-square property holds:
for any closed cycle $\gamma$ of length 4 in the 1-skeleton of $K$
there is at least one diagonal of $\gamma$ (i.e. an edge connecting
some pair of opposite vertices of $\gamma$) in $K$.
Recall that a right-angled Coxeter group $W$ is word-hyperbolic if
and only if  its nerve $L$
satisfies the flag-no-square condition
(see Corollary 12.6.3 in [Da], as applied to the right angled case).

We will first show an example of a 3-sphere
with holes  which is flag-no-square and which satisfies all the assumptions
of Theorem 1.3. 
More precisely, by Lemma 4.8
below, the 580-disk described in Definition 4.7 is such an example.
In view of Theorem 1.3, this leads to an example
of a right-angled hyperbolic Coxeter group whose Gromov
boundary is homeomorphic to the 2-dimensional Sierpi\'nski compactum.
On the other hand, we will show that for $n\ge3$ there are no examples of 
$(n+1)$-spheres with holes which satisfy the flag-no-square condition,
and to which Theorem 1.3 applies (see Lemma 4.9 below).

Consider the 600-cell, the convex regular 4--polytope with Schl\"{a}fli symbol $\{3,3,5\}$ (see e.g. [C]).
Denote by $X_{600}$ the boundary of the
600--cell, a~3-dimensional
simplicial polyhedron homeomorphic to the 3-dimensional sphere.
Recall that, when viewed as a simplicial complex, $X_{600}$
consists of 600 3-simplices and has 120 vertices. Its vertex links
are icosahedra and edge links are pentagons. We will only
exploit the combinatorial simplicial structure of $X_{600}$.

\medskip\noindent
{\bf 4.7 Definition.} The {\it 580-disk} is the simplicial complex
obtained from $X_{600}$ by deleting open star of any vertex.
We denote it by $X_{580}$.

\medskip
Note that $X_{580}$ is a PL triangulation of the 3-dimensional disk $B^3$ consisting
of 580 3-simplices, and its boundary $\partial X_{580}$, which corresponds to the link of the deleted vertex, is an icosahedral triangulation of the 2-sphere. 

\medskip\noindent
{\bf 4.8 Lemma.}
{\it The 580-disk $X_{580}$ is a flag-no-square simplicial complex, 
and the boundary $\partial X_{580}$
is its 3-convex subcomplex.}

\medskip\noindent
{\bf Proof:}
To prove the lemma, we need the following fact, whose proof can be
found in [PS] 
(see Lemma 2.1 in that paper).

\medskip\noindent
{\bf Lemma 1.}
{\it $X_{600}$ satisfies the flag-no-square property.}

\medskip
Next two claims exhibit basic properties of the 580-disk,
resulting from Lemma 1.

\medskip\noindent
{\bf Claim 1.}
{\it The 580-disk $X_{580}$ satisfies the flag-no-square property.}

\medskip
To get Claim 1,
note that $X_{580}$ is a full subcomplex in $X_{600}$.
This follows from a more general easy observation that, given any
simplicial complex $K$ and its any vertex $v$, the subcomplex of $K$ obtained
by deleting the open star of $v$ is full. Since the flag-no-square property is
inherited by full subcomplexes, the claim follows by referring to
Lemma 1.

\medskip\noindent
{\bf Claim 2.}
{\it The boundary $\partial X_{580}$ is a full subcomplex of the 580-disk $X_{580}$.}

\medskip
To get Claim 2, denote by $D$ the closed star of this vertex of 
$X_{600}$ which does not belong to $X_{580}$. Since closed star
of any vertex in a flag simplicial complex is a full subcomplex,
and since $X_{600}$ is flag (see Lemma 1), we get that $D$ is a full
subcomplex of $X_{600}$. 
Since 
$\partial X_{580}=D\cap X_{580}$, the claim is a consequence
of the following general fact: 
if $A$ is a full subcomplex of
a simplicial complex $K$ then for any subcomplex $B$ of $K$ the intersection
$A\cap B$ is a full subcomplex of $B$.

\medskip
In view of Claims 1 and 2, it remains to show that, given any two
vertices $x,z$ of $\partial X_{580}$ at combinatorial distance 2
in the 1-skeleton of $X_{580}$, and any polygonal path $(x,y,z)$
connecting them in $X_{580}$, we have $y\in \partial X_{580}$.
View $X_{580}$ as a subcomplex in $X_{600}$, and let $v$ be the vertex
of $X_{600}$ not contained in $X_{580}$. Consider the 4-cycle
$(x,y,z,v,x)$ in $X_{600}$. Since the combinatorial distance
between $x$ and $z$ in the 1-skeleton of $X_{600}$ is easily seen to be 2
(e.g. since $X_{580}$ is full in $X_{600}$),
it follows from Lemma 1 that $y$ is adjacent to $v$ in $X_{600}$.
But this means that $y\in\partial X_{580}$, which completes the proof of Lemma 4.8.

\medskip
Our final observation is the following lemma, which implies that
using Theorem 1.3, one cannot construct an example of a right-angled word-hyperbolic Coxeter group whose Gromov boundary is an $n$-dimensional
Sierpi\'nski compactum, for any $n\ge3$. 

\medskip\noindent
{\bf 4.9 Lemma.}
{\it If $n\ge3$, then there is no $(n+1)$-sphere with holes $L$ 
satisfying conditions (1) and (2) of Theorem 1.3 and
the flag-no-square
condition.}

\medskip\noindent
{\bf Proof:}
Let $L$ be an $(n+1)$-sphere with holes satisfying the assumtions (1) and (2)
of Theorem 1.3 and the flag-no-square
condition. For any peripheral cycle $C$ of $L$, attach to $L$ the simplicial
cone over $C$, thus getting a PL triangulation $L'$ of the $(n+1)$-sphere.
We will show that $L'$ satisfies the flag-no-square condition.
 
Flagness of $L'$ easily follows from flagness of $L$ and the fact
that the peripheral cycles are full subcomplexes of $L$.
It remains to show that $L'$ satisfies the no-empty-square condition.
More precisely,
let $\gamma=(w,x,y,z,w)$ be a 4-cycle in the 1-skeleton of $L'$.
We need to show that $\gamma$ has at least one of the diagonals in $L'$.

If $\gamma$ is contained in $L$, existence of a diagonal follows from
the assumption that $L$ is flag-no-square. Otherwise, $\gamma$ has
a vertex, say $w$, which is the cone vertex of one of the cones,
say the cone over a peripheral cycle $C$, attached to $L$.
We then have that $x$ and $z$ necessarily belong to $C$.
If $y$ is also a vertex of $C$, the edge $(w,y)$ is a diagonal
of $\gamma$. If $x$ is adjacent to $z$ in $C$, the edge $(x,z)$
is a diagonal of $\gamma$. 
It remains to consider the last case when $x$ is not adjacent to $z$ in $C$
and $y$ is not contained in $C$.

The last case mentioned above has two subcases.
In the first subcase, suppose that $y$ is a vertex of $L$.
Then existence of the path $(x,y,z)$ contradicts 3-convexity of $C$ in $L$,
and hence this subcase is not possible.
In the second and last subcase, $y$ has to be the cone vertex of the cone
over a peripheral cycle $C'$ of $L$ distinct from $C$.
But then the non-adjacent vertices $x,z$ belong to the intersection
$C\cap C'$, which contradicts the assumption that such intersection
is either empty or a single simplex of $L$. So this subcase is also not possible,
and hence $L'$ satisfies the flag-no-square condition.

The lemma is now a consequence of nonexistence of flag-no-square
triangulations of spheres $S^{n+1}$ for $n\ge3$
(cf. Corollary 5.7 in [PS]).

\medskip
Lemma 4.9, and the comment right before its statement, suggest the
following.

\medskip\noindent
{\bf 4.10 Question.}
Given an integer $n\ge3$, does there exist a right-angled Coxeter group $W$
which is word-hyperbolic and whose Gromov boundary $\partial W$
is an $n$-dimensional Sierpi\'nski compactum?


\bigskip
\noindent
{\bf 5. Local 3-convexity in \hbox{CAT}(0) cubical complexes.} 

\medskip
In this section we use the concept of 3-convexity 
to present some strengthened notion of convexity (called {\it local 3-convexity})
in CAT$(0)$ cubical complexes. We then derive, as Proposition 5.3 below, 
some rather strong global geometric 
property of locally 3-convex subcomplexes. This property is
essential for our arguments in Section 6.

\medskip\noindent
{\bf 5.1 Definition} (local 3-convexity).
A subcomplex $Y$ of a cubical complex $X$ is {\it locally 3-convex}
if for each vertex $v\in Y$ the link $\hbox{Lk}(v,Y)$ is a 3-convex
subcomplex in the corresponding link $\hbox{Lk}(v,X)$.

\medskip
Before stating the main result of this section, concerning some global behaviour
of locally 3-convex subcomplexes, we make few preparations.

To each complete convex subspace $Y$ of a CAT(0) space $X$
(in particular, to any convex subcomplex $Y$ of a CAT(0) cubical
complex $X$) there is associated a  map $\pi:X\to Y$, 
called {\it projection to $Y$}, given by $\pi(x)=y$, where $y$ is the unique
point of $Y$ closest to $x$ (cf. Chapter II.2 in [BH]). The following 
observation will be used later in our arguments.

\medskip\noindent
{\bf 5.2 Lemma.}
{\it Let $X$ be a CAT(0) cubical complex and $Y$ its convex subcomplex.
Then for any vertex $v$ of $X$ its projection to $Y$ is also a vertex.}

\medskip\noindent
{\bf Proof:} Supose that $\pi(x)=y$, where $y\in Y$ is not a vertex.
We need to show that $x$ is not a vertex.
Let $\sigma$ be the cubical face of $Y$ such that $y$ belongs to the
interior of $\sigma$. Let $[x,y]$ denote the geodesic in $X$ from $x$ to $y$.
Since $y$ is the point of $Y$ closest to $x$, there is a face $\tau$ in $X$
containing $\sigma$ and such that the final part of $[x,y]$ is contained in $\tau$ and orthogonal to $\sigma$. 
In particular, there is a pair of opposite codimension 1 faces in $\tau$
such that the part of $[x,y]$ contained in $\tau$ is parallel to those
two faces, and disjoint with them both.
Then, using induction, we observe that $[x,y]$ satisfies the same property
with respect to any cubical face $\tau'$ of $X$ that it crosses
(here we use the fact that $[x,y]$ is the concatenation of finitely many segments, each of which is contained in a single cubical face,
cf. Corollary 7.29 on p. 110 in [BH]). In particular, it follows that
$[x,y]$ omits
all vertices of any cubical face of $X$ that it crosses, and hence $x$ is not
a vertex of $X$. This completes the proof.

\medskip
We denote by $\hbox{st}_c(v,X)$, and call the {\it cubical star} 
of a vertex $v$ in
a CAT(0) cubical complex $X$, the subcomplex 
of $X$ equal to the union of all cubical faces of $X$ containg $v$.

\medskip\noindent
{\bf 5.3 Proposition.}
{\it Let $X$ be a CAT$(0)$ cubical complex, and let $Y$ be its convex
subcomplex which is locally 3-convex. Let $v$ be any vertex of $X$,
and denote by $u$ the vertex of $Y$ equal to the projection of $v$ to $Y$.
Let $\rho$ be any geodesic ray in $X$ started at $v$ and representing
a point $\xi\in\partial_\infty Y\subset\partial_\infty X$. Then}
\item{(a)} {\it $\rho$ intersects the subcomplex $\hbox{st}_c(u,Y)$;}
\item{(b)} {\it except for a bounded initial subsegment, $\rho$ 
is contained in $Y$.}

\medskip\noindent
{\bf Proof:} The proof is split into several parts, and it occupies
the rest of this section.

\medskip\noindent
{\it Thickened hyperplanes and extended half-spaces.}

We refer to hyperplanes and half-spaces in CAT(0) cubical complexes,
as described e.g. in [W], Sections 2.4, 3.2 and 3.3.
Given any oriented edge $e=(p,q)$ of $X$, denote by $J_e$
the hyperplane in $X$ dual to $e$, i.e. the hyperplane orthogonal to $e$
and passing through the midpoint of $e$. Denote
by $H_e$ this closed half-space
for $J_e$ which contains $p$. Recall that both $J_e$ and $H_e$ are
convex subspaces in $X$, and that the hyperplane $J_e$ 
is canonically a CAT(0) cubical complex. 

Denote by $J_e^*$, and call the {\it thickened hyperplane} in $X$ dual
to $e$, the subcomplex of $X$ being the union of all cubical faces of $X$
intersected by $J_e$. Note that, as a cubical complex, $J_e^*$ is 
isomorphic to the product $J_e\times I$, where $I$ is the 1-cube.
Moreover, given any vertex $o$ in $J_e^*$, there is a unique edge $\varepsilon$ in $X$ containing $o$ and intesecting $J_e$; denoting
by $t_\varepsilon$ the vertex in the link $\hbox{Lk}(o,X)$ induced by 
$\varepsilon$, we have $\hbox{Lk}(o,J_e^*)=
\hbox{st}(t_\varepsilon,\hbox{Lk}(o,X))$.
In particular, $\hbox{Lk}(o,J_e^*)$ is a full subcomplex of $\hbox{Lk}(o,X)$,
and hence any thickened hyperplane $J_e^*$ is a convex subcomplex of $X$.

Denote by $H_e^*$, and call the {\it extended half-space} of $X$
induced by $e$, the subcomplex equal (as a set) to the union $H_e\cup J_e^*$.
Since for any vertex $s\in H_e^*$ either 
$\hbox{Lk}(s,H_e^*)=\hbox{Lk}(s,X)$ or 
$\hbox{Lk}(s,H_e^*)=\hbox{Lk}(s,J_e^*)$,
we get that any extended half-space $H_e^*$ is a convex subcomplex of $X$.

Compare [W], Section 3.3, where extended hyperplanes are called hyperplane carriers, and extended half-spaces are called halfspace carriers.

\medskip\noindent
{\it The subcomplex $P$ and its basic properties.}

Given $X$, $Y$, $v$ and $u$ as in the proposition,
denote by ${\cal E}_u^Y$ the set of all oriented edges in $Y$
with the initial vertex $u$, and put 
$$
P:=\bigcap_{e\in{\cal E}_u^Y}H_e^*.
$$
Obviously, $P$ is a convex subcomplex of $X$.

\medskip\noindent
{\bf Claim 1.} $P\cap Y=\hbox{st}_c(u,Y)$.

\medskip
To prove Claim 1, note that for each $e\in{\cal E}_u^Y$,
if we denote by $J_{e,Y}$ the hyperplane in $Y$ dual to $e$, we have
$J_{e,Y}=J_e\cap Y$. Similarly, if $H^*_{e,Y}$ denotes the extended
half-space in $Y$ induced by $e$, we have $H^*_{e,Y}=H_e^*\cap Y$.
Consequently, we have
$$
P\cap Y=\bigg(\bigcap_{e\in{\cal E}_u^Y}H_e^*\bigg)\cap Y=
\bigcap_{e\in{\cal E}_u^Y}(H_e\cap Y)=
\bigcap_{e\in{\cal E}_u^Y}H^*_{e,Y},
$$
and the latter intersection obviously coincides with the cubical star
$\hbox{st}_c(u,Y)$, thus justifying the claim.

\medskip\noindent
{\bf Claim 2.} $v\in P$.

\medskip
Since $H_e\subset H_e^*$,
to prove Claim 2, it is sufficient to show that for each
$e\in{\cal E}_u^Y$ we have $v\in H_e$. To do this, suppose it is not true,
i.e. $v\notin H_e$ for some $e\in{\cal E}_u^Y$.
Since $u\in H_e$, the geodesic $[v,u]$ crosses the
hyperplane $J_e$. Let $z$ be the point in the intersection $[v,u]\cap J_e$
closest to $u$. Consider the geodesic triangle $T=\Delta(u,m,z)$, where $m$
is the midpoint of the edge $e$. Since $z\in J_e$, the angle of $T$ at $m$
equals $\pi/2$. Since the side $[u,m]$ of $T$ is contained in $Y$, and since
$u$ is the projection of $z$ to $Y$, the angle of $T$ at $u$ is $\ge\pi/2$.
Finally, the angle of $T$ at $z$ is $>0$.
This contradicts the fact that the sum of angles in any geodesic triangle
in a CAT(0) space is $\le\pi$, thus completing the proof of Claim 2.

\medskip\noindent
{\it Convexity of $Y\cup P$.}

To prove that the subcomplex $Y\cup P$ is convex, 
observe first that it is obviously connected.
It is then sufficient to show that for any vertex 
$z\in Y\cup P$ the link $\hbox{Lk}(z,Y\cup P)$ is a full subcomplex
in the link $\hbox{Lk}(z,X)$. Obviously, this is true for any vertex 
$z\in Y\cup P$ disjoint with the intersection $P\cap Y$,
by convexity of both $P$ and $Y$. It remains to check this property
for vertices $z\in P\cap Y$, i.e. for vertices $z\in\hbox{st}_c(u,Y)$,
see Claim 1 above.

For $z=u$ we have $\hbox{Lk}(u,P\cup Y)=\hbox{Lk}(u,X)$,
because $\hbox{st}_c(u,X)\subset P\cup Y$, and hence the above mentioned
property follows
trivially. Suppose $z\in P\cap Y=\hbox{st}_c(u,Y)$, 
and $z\ne u$. Let $\sigma$ be the smallest
cubical face of $X$ containing $u$ and $z$.
Let $e_1,\dots,e_m$ be the edges adjacent to $u$ in $\sigma$,
oriented so that $u$ is their initial vertex. Obviously, we have 
$e_1,\dots,e_m\in{\cal E}_u^Y$. Let $d_1,\dots,d_m$ be the edges
of $\sigma$ adjacent to $z$ and parallel to $e_1,\dots,e_m$ respectively.
Equip them with the orientations consistent with those of $e_1,\dots,e_m$,
i.e. such that $z$ is their terminal vertex.
Denote by $t_1,\dots,t_m$ the vertices in $\hbox{Lk}(z,X)$ induced by
the edges $d_1,\dots,d_m$ respectively. Denote also by $\tau$ the simplex
of $\hbox{Lk}(z,X)$ spanned by the vertices $t_1,\dots,t_m$
(i.e. induced by the face $\sigma$). We then have
$$
\hbox{Lk}(z,P)=\hbox{Lk}\bigg(z,\bigcap_{e\in{\cal E}_u^Y}H_e^*\bigg)=
\bigcap_{e\in{\cal E}_u^Y}\hbox{Lk}(z,H_e^*)=
\hbox{Lk}(z,H_{e_1}^*)\cap\dots\cap\hbox{Lk}(z,H_{e_m}^*),
$$
where the last equality follows by an easy observation that for 
$e\in{\cal E}_u^Y\setminus\{ e_1,\dots,e_m \}$ we have
$\hbox{Lk}(z,H_e^*)=\hbox{Lk}(z,X)$.
As a consequence, we get
$$
\hbox{Lk}(z,P\cup Y)=\hbox{Lk}(z,Y)\cup\hbox{Lk}(z,P)=
\hbox{Lk}(z,Y)\cup[\hbox{Lk}(z,H_{e_1}^*)\cap\dots\cap\hbox{Lk}(z,H_{e_m}^*)]=
$$
$$
=\hbox{Lk}(z,Y)\cup[\hbox{Lk}(z,H_{d_1}^*)\cap\dots\cap\hbox{Lk}(z,H_{d_m}^*)]=
$$
$$
=\hbox{Lk}(z,Y)\cup[\hbox{Lk}(z,J_{d_1}^*)\cap\dots\cap\hbox{Lk}(z,J_{d_m}^*)]=
$$
$$
=\hbox{Lk}(z,Y)\cup[\hbox{st}(t_1,\hbox{Lk}(z,X))\cap\dots\cap 
\hbox{st}(t_m,\hbox{Lk}(z,X))]=
\hbox{Lk}(z,Y)\cup\hbox{st}(\tau,\hbox{Lk}(z,X)).
$$
Since $\hbox{Lk}(z,Y)$ is 3-convex in $\hbox{Lk}(z,X)$,
and since $\tau\subset\hbox{Lk}(z,Y)$, it follows from the above equality
and from the second assertion 
of Lemma 4.2 that $\hbox{Lk}(z,P\cup Y)$ is a full subcomplex of
$\hbox{Lk}(z,X)$. This completes the proof of
convexity of the subcomplex $P\cup Y\subset X$.

\medskip\noindent
{\it Proof of assertions (a) and (b) of the proposition.} 

Let $\rho$ be a geodesic ray as in the proposition.  
If $v\in Y$,  assertion (a) is obvious since $u=v$, and assertion (b) follows
from convexity of $Y$. So assume that $v\notin Y$.
By the fact that $Y\cup P$ is convex, it makes sense to speak of
the subset $\partial_\infty(Y\cup P)\subset\partial_\infty X$. Since 
$\partial_\infty Y\subset\partial_\infty(Y\cup P)$, we have $\xi\in\partial_\infty(Y\cup P)$, and since $v\in P\subset Y\cup P$,
we get that $\rho$ is contained in $Y\cup P$.
To prove (a), suppose on the contrary that $\rho$ does not intersect 
$\hbox{st}_c(u,Y)=P\cap Y$. 
Since $P\setminus(P\cap Y)$ is open and closed in $(P\cup Y)\setminus(P\cap Y)$,
and since $v\in P\setminus(P\cap Y)$, the ray $\rho$ must be contained in
$P\setminus(P\cap Y)$. 
Since $P\cap Y$ is bounded, it follows that 
$\hbox{d}_{Y\cup P}(\rho(t),Y)=\hbox{d}_{X}(\rho(t),Y\cap P)\to\infty$ as $t\to\infty$.
By convexity of $Y\cup P$, the distances in $Y\cup P$ are the same
as those in $X$, and hence 
$\hbox{d}_{X}(\rho(t),Y)\to\infty$ as $t\to\infty$,
which contradicts the assumption that 
$\lim_{t\to\infty}\rho(t)=\xi\in\partial_\infty Y$.
Thus assertion (a) follows.
A similar argument shows that, after leaving $P\cap Y$, the geodesic ray
$\rho$ stays in $Y$, which yields assertion (b).


\bigskip
\noindent
{\bf 6. Completing the proof of Theorem 1.3.}
\medskip

To prove Theorem 1.3, we consider 
(under notation and assumptions of this theorem)
an embedding of the visual boundary
$\partial_\infty(W,S)$ in the sphere $S^{n+1}$ as described in Section 2.
We need to verify that, under this embedding, 
$\partial_\infty(W,S)$ satisfies conditions (1)--(4) of Definition 1.1.
Note that conditions (1) and (2) have been already verified in Section 2
(see Proposition 2.1(2) and Corollary 2.6). It remains to verify
conditions (3) and (4), which we do below, 
in Propositions 6.1 and 6.4.


\medskip
\noindent
{\it Closures of complementary components are pairwise disjoint.}

\smallskip
This subsection is devoted to the proof of the following.

\medskip\noindent
{\bf 6.1 Proposition.}
{\it Let $(W,S)$ be a right-angled Coxeter system whose nerve $L$ is an
$(n+1)$-sphere with holes, and suppose that
the following two conditions hold:}
\item{(1)} {\it for any two distinct connected components of $S^{n+1}\setminus L$,
their closures in $S^{n+1}$ are either disjoint, or intersect at
a single simplex of $L$;}
\item{(2)} {\it the boundary $\hbox{\rm bd}(\Omega)$ of any component $\Omega$ of $S^{n+1}\setminus L$ is a 3-convex 
subcomplex of $L$.}

\noindent
{\it Consider an embedding of the boundary $\partial_\infty(S,W)$ in $S^{n+1}$,
as described in Section 2. Then the closures of any two distinct connected components of 
$S^{n+1}\setminus\partial_\infty(W,S)$ are disjoint.}

\medskip
We use the setting and the notation fixed at the
beginning of Section 2. In particular, $\partial_\infty(W,S)=\partial_\infty\Sigma$ is realized as a subspace of
$S^{n+1}=\partial_\infty\Sigma'=\partial_\infty(W',S')$,
for appropriate choice of the nerve $L'\supset L$. By Lemma 2.5,
connected components of $S^{n+1}\setminus\partial_\infty(W,S)=
\partial_\infty\Sigma'\setminus\partial_\infty\Sigma$
are parametrized by the pairs $(C,\bar w):C\in\hbox{Per}(L,L'), \bar w\in W/W_C$ and have form
$$
U=U_{C,\bar w}=\partial_\infty\Sigma'\setminus\partial_\infty(\bar w\cdot H_C^+)
=\partial_\infty(\bar w\cdot H_C^-)\setminus\partial_\infty(\bar w\cdot\Sigma_C).
$$
The topological boundary $\hbox{bd}(U)$ coincides then with the set
$\partial_\infty(\bar w\cdot\Sigma_C)$. Thus, to prove Proposition 6.1, it is
sufficient to show the following.

\medskip\noindent
{\bf 6.2 Lemma.}
{\it Under assumptions of Proposition 6.1,
if $(C_1,\bar w_1), (C_2,\bar w_2)$ are distinct pairs such that 
$C_i\in\hbox{Per(L,L')}$ and $\bar w_i\in W/W_{C_i}$ then the sets
$\partial_\infty(\bar w_1\cdot\Sigma_{C_1})$ and $\partial_\infty(\bar w_2\cdot\Sigma_{C_2})$ are disjoint (as subsets of $\partial_\infty\Sigma$).}

\medskip\noindent
{\bf Proof:}
We will need the following.

\medskip\noindent
{\bf Claim.}
{\it If the intersection $\bar w_1\cdot\Sigma_{C_1}\cap\bar w_2\cdot\Sigma_{C_2}$ is not empty
then it consists of a single cubical cell of the complex $\Sigma$.}

\medskip
To prove the claim,
we refer here to the combinatorial description of the Coxeter-Davis 
complex $\Sigma$, as a cubical complex,
given in [Da], Chapter 7, Section 7.3.
In this description, the cubical cells of $\Sigma$ are in 1-1 
correspondence with the cosets of finite special subgroups of $W$,
and in particular the vertices of $\Sigma$ are in 1-1 correspondence
with $W$.
For any $C\in\hbox{Per}(L,L')$ and any $\bar w=w\cdot W_C\in W/W_C$, the vertex set of the subcomplex $\bar w\cdot\Sigma_C$
can be identified with the subset $w\cdot W_C\subset W$.
Obviously, if $C_1=C_2=C$ and $\bar w_1=w_1\cdot W_C\ne\bar w_2=w_2\cdot W_C$ then the subsets  $w_1\cdot W_{C_1}$ and
$w_2\cdot W_{C_2}$ are disjoint. It follows that under assumption
of the claim we have $C_1\ne C_2$.
Now, since $\bar w_1\cdot\Sigma_{C_1}\cap\bar w_2\cdot\Sigma_{C_2}\ne\emptyset$,
there is a vertex $w\in W$ which belongs to this intersection.
We then have $\bar w_1\cdot\Sigma_{C_1}=w\cdot\Sigma_{C_1}$ and
$\bar w_2\cdot\Sigma_{C_2}=w\cdot\Sigma_{C_2}$.
Consequently, by invariance under $W$, to prove the claim it is sufficient to show that the intersection $\Sigma_{C_1}\cap\Sigma_{C_2}$
is a single cell.
The vertex sets of the subcomplexes $\Sigma_{C_1}$ and $\Sigma_{C_2}$ correspond to
the subgroups $W_{C_1},W_{C_2}<W$, respectively. Obviously,
we have $W_{C_1}\cap W_{C_2}=W_{C_1\cap C_2}$, where 
$W_\emptyset$ denotes the trivial subgroup. By the assumption (1)
of Proposition 6.1, $C_1\cap C_2$ is either empty, or a single simplex
of $L$. The special subgroup $W_{C_1\cap C_2}$ is then finite.
As a consequence, the intersection $\Sigma_{C_1}\cap\Sigma_{C_2}$
is exactly the cell of $\Sigma$ corresponding to the subgroup
$W_{C_1\cap C_2}$ viewed as a coset of a finite special subgroup in $W$,
hence the claim.

\medskip
To conclude the proof of Lemma 6.2,
consider any two points $\xi_1,\xi_2$ belonging to the sets 
$\partial_\infty(\bar w_1\cdot\Sigma_{C_1})$ and $\partial_\infty(\bar w_2\cdot\Sigma_{C_2})$, respectively. Fix any vertex $v$ of $\Sigma$,
and denote by $\rho_1,\rho_2$ the geodesic rays in $\Sigma$
started at $v$ and representing the points $\xi_1,\xi_2$ respectively.
Note that for any vertex $u$ of any subcomplex $\bar w_i\cdot\Sigma_{C_i}$
the pair of vertex links 
$(\hbox{Lk}(u,\Sigma),\hbox{Lk}(u,\bar w_i\cdot\Sigma_{C_i})$
is isomorphic (as a pair of simplicial complexes) to the pair $(L,C_i)$.
By assumption (2) of Proposition 6.1, it follows that each
$\bar w_i\cdot\Sigma_{C_i}$ is a convex and locally 3-convex subcomplex
of $\Sigma$. Applying Proposition 5.3(b), we get that each of the rays $\rho_i$, except for a finite subsegment, is contained in the corresponding
subcomplex $\bar w_i\cdot\Sigma_{C_i}$.
In view of the claim, this implies that $\rho_1\ne\rho_2$,
and consequently $\xi_1\ne\xi_2$, hence the lemma.

\medskip
By what was said above, this also concludes the proof of Proposition 6.1.

\medskip
\noindent
{\it Natural metrics on the visual boundary of a $\hbox{CAT}(0)$
space.}

\smallskip
In this subsection we describe a family of natural metrics
(compatible with the cone topology) on the visual boundary $\partial_\infty X$
of a complete CAT(0) space $X$. 
(A metric of this form will be used in the argument in the next subsection.)
Recall that in a complete CAT(0) space
any point $x_0\in X$ can be connected to any $\xi\in\partial_\infty X$
with the unique geodesic ray (see [BH], Chapter II, Proposition 8.19).
Recall also that for distinct geodesic rays $r,r'$ based at a common point of $X$
the function $t\to d_X(r(t),r'(t))$ is convex (see [BH], Chapter II, Proposition 2.2),
and hence for any real $A>0$ there is exactly one $t$ with $d_X(r(t),r'(t))=A$.
The following fact was observed by Damian Osajda, see Proposition 9.6(1)
in [OS].

\medskip\noindent
{\bf 6.3 Lemma.}
{\it Let $X$ be a complete CAT(0) space, $x_0\in X$ a point,
and $A>0$ a real number. For $\xi,\eta\in\partial_\infty X, \xi\ne\eta$, let
$r_\xi,r_\eta$ be the unique geodesic rays from $x_0$ to $\xi$ and $\eta$, respectively.
Put $d_A(\xi,\eta):=t^{-1}$, where $t$ is the unique number with $d_X(r_\xi(t),r_\eta(t))=A$. Then $d_A$ is a metric on $\partial_\infty X$ compatible with the cone topology.}

\bigskip
\noindent
{\it The family of complementary components is null.}

\smallskip
The last missing element in the proof of Theorem 1.3 is the following.

\medskip\noindent
{\bf 6.4 Proposition.}
{\it Let $(W,S)$ be a right-angled Coxeter system whose nerve $L$ is an
$(n+1)$-sphere with holes, and suppose that
the following two conditions hold:}
\item{(1)} {\it for any two distinct connected components of $S^{n+1}\setminus L$,
their closures in $S^{n+1}$ are either disjoint, or intersect at
a single simplex of $L$;}
\item{(2)} {\it the boundary $\hbox{\rm bd}(\Omega)$ of any component $\Omega$ of $S^{n+1}\setminus L$ is a 3-convex 
subcomplex of $L$.}

\noindent
{\it Consider an embedding of the boundary $\partial_\infty(S,W)$ in $S^{n+1}$,
as described in Section 2. Then the family of connected components of 
$S^{n+1}\setminus\partial_\infty(W,S)$ is null.}

\medskip\noindent
{\bf Proof:} We come back to the notation introduced at the beginning 
of Section 2. Note that, in view of Corollary 2.6, the family of closures of the connected
components of $S^{n+1}\setminus\partial_\infty(W,S)$ coincides with
the family
$$
\partial_\infty(\bar w\cdot H_C^-): C\in \hbox{Per}(L,L'), \bar w\in W/W_C.
$$
Obviously, it is then sufficient to show that this last family is null.

Order the family of all half-spaces in $\Sigma'$ of form $\bar w\cdot H_C^-$
as above into a sequence $(H_k)_{k\ge1}$, and recall that these half-spaces
are all convex.
Fix any vertex $v_0\in\Sigma'$, and for each $k\ge1$ denote by
$u_k$ the projection of $v_0$ to $H_k$ (which is a vertex, by Lemma 5.2).

\medskip\noindent
{\bf Claim 1.}
{\it Put $t_k=d_{\Sigma'}(v_0,u_k)$. Then $\lim_{k\to\infty}t_k=\infty$.}

\medskip
To prove Claim 1, note that each half-space $H_k$ is the union of a collection
of $(n+2)$-cubes of $\Sigma'$, and that two different half-spaces $H_k$
have no such $(n+2)$-cube in common. The claim follows then from
uniform finiteness of $\Sigma'$.

\medskip\noindent
{\bf Claim 2.}
{\it Each half-space $H_k$ is a convex and
locally 3-convex subcomplex of $\Sigma'$.}

\medskip
To prove Claim 2, we only need to check that $H_k$ is locally 3-convex
in $\Sigma'$.
Assume without loss of generality that 
$H_k=\bar w\cdot H_C^-$, and put $\Sigma_k=\bar w\cdot\Sigma_C$.
Obviously, if $v$ is a vertex of $H_k$ not belonging to $\Sigma_k$
then the link $\hbox{Lk}(v,H_k)$ coincides with the link 
$\hbox{Lk}(v,\Sigma')$. If $v\in\Sigma_k$
then the pair $(\hbox{Lk}(v,\Sigma'), \hbox{Lk}(v,H_k))$
is isomorphic to the pair $(L',\overline\Omega)$,
where $\overline\Omega$ is the subcomplex of $L'$
equal to the closure of the connected component of $L'\setminus L$
bounded by $C$. Thus, in order to prove the claim,
we need to show that the subcomplex
$\overline\Omega$ is 3-convex in $L'$.

To prove the latter statement, note first that since $C$ is a full subcomplex
of $L'$, the same is true for $\overline\Omega$. To get 3-convexity
of $\overline\Omega$, consider any combinatorial geodesic $(x,y,z)$
in the 1-skeleton of $L'$, with $x$ and $z$ belonging to $\overline\Omega$.
We need to show that $y$ belongs to $\overline\Omega$.
This is obviously true if some of the vertices $x,z$ does not lie
on the boundary cycle $C$. So we may assume that both $x$ and $z$
are in $C$. If $y$ is a vertex of $L$, it must belong to $C$, and hence to
 $\overline\Omega$, by 3-convexity of $C$ in $L$. If $y$ is not in $L$,
it belongs to some component $\Omega'$ of the complement $L'\setminus L$.
It follows that both $x$ and $z$ belong then to the closure
$\overline\Omega'$.  Since $x$ and $z$ are not adjacent,
condition (1) of the proposition implies that $\Omega'=\Omega$,
and hence again $y\in\overline\Omega$.
This completes the proof of Claim 2.

\medskip
To conclude the proof of Proposition 6.4,
we now pass to estimating the diameter of the subset $\partial_\infty H_k$
with respect to some metric $d_A$ in $\partial_\infty\Sigma'=S^{n+1}$.
More precisely, 
put $A=2\sqrt{n+2}$, and consider the metric $d_A$ in $\partial_\infty\Sigma'$ as described in Lemma 6.3, for the base point $v_0$.
Let $\xi,\eta$ be any points of $\partial_\infty H_k$, and let 
$r_\xi,r_\eta$ be the geodesic rays in $\Sigma'$, started at $v_0$, and
representing points $\xi,\eta$ respectively.
In view of Claim 2, Proposition 5.3(a) implies that $r_\xi,r_\eta$ 
pass through the cubical star  $\hbox{st}_c(u_k,H_k)$.
Since obviously for each $k\ge1$ the diameter
of the cubical star $\hbox{st}_c(u_k,H_k)$ is bounded above by $A$,
and since the distance from $v_0$ to $H_k$ is equal to $t_k$,
we get that $d_{\Sigma'}(r_\xi(t_k),r_\eta(t_k))\le A$,
and hence $d_A(\xi,\eta)\le1/t_k$.
This means that for the metric $d_A$ we have 
$\hbox{diam}(H_k)\le1/t_k$.
The proposition follows then by Claim 1.


\bigskip\noindent
{\bf References}

\medskip

\item{[Ar]} S. Armentrout, {\it Cellular decompositions of 3-manifolds
that yield 3-manifolds}, Mem. Amer. Math. Soc. 107 (1971).

\item{[Be]} M. Bestvina, {\it Local homology properties of boundaries of groups,}
Michigan Math, J. 43 (1996), 123--139.

\item{[BeM]} M. Bestvina, G. Mess, {\it The boundary of negatively curved
groups}, J. Amer. Math. Soc. 4 (1991), 469--481.

\item{[BH]} M. Bridson, A. Haefliger, 
Metric Spaces of Non-Positive
Curvature, Grundlehren der mathematischen Wissenschaften 319, Springer,
1999.

\item{[Br]} M. Brown, {\it Some applications of an approximation
theorem for inverse limits,} Proc. Amer. Math. Soc. 11 (1960), 478--481. 

\item{[Ca]} J. Cannon, {\it A positional characterization of the
$(n-1)$-dimensional Sierpi\'nski curve in $S^n$,} Fund. Math. 79 (1973),
107--112.

\item{[C]} H.S.M. Coxeter, {Regular polytopes}, Methuen and Co., London, 1948.

\item{[DV]} R. Daverman, G. Venema, Embeddings of Manifolds,
Graduate Studies in Mathematics, vol. 106, American Mathematical Society,
Providence, Rhode Island, 2009.

\item{[Da]} M. Davis, {The geometry and topology 
of Coxeter groups}, {London Mathematical Society Monographs Series},
{vol. 32}, {Princeton University Press}, {Princeton}, {2008}.

\item{[DJ]} M. Davis, T. Januszkiewicz, {\it Hyperbolization of polyhedra},
J. Diff. Geom. 34 (1991), 347--388.

\item{[Dr]} A. Dranishnikov, {\it On boundaries of hyperbolic Coxeter groups}, Topology and its Applications 110 (2001), 29--38. 

\item{[He]} J. P. Henderson, {\it Cellularity in polyhedra,}
Topology and its Applications 12 (1981), 267--282.

\item{[Mo]} R. L. Moore, {\it Concerning upper semi-continuous
collections of continua}, Trans. Amer. Math. Soc. 27 (1925), 416--428.

\item{[OS]} D. Osajda, J. \'Swi\c atkowski, {\it  On asymptotically hereditarily
aspherical groups}, Proc. London Math. Soc. 111 (2015),  93--126.

\item{[PS]} P. Przytycki, J. \'Swi\c atkowski, 
{\it Flag-no-square triangulations and Gromov boundaries
in dimension 3}, Groups, Geometry \& Dynamics  3 (2009), 453--468.

\item{[Q]} F. Quinn, {\it Ends of maps. III: dimensions 4 and 5,}
J. Diff. Geom. 17 (1982), 503--521.

\item{[Si]} L. C. Siebenmann, {\it Approximating cellular maps by homeomorphisms,} Topology 11 (1972), 271--294.

\item{[Sw1]} J. \'Swi\c atkowski, {\it Trees of manifolds as boundaries
of spaces and groups,} preprint 2013,  arXiv:1304.5067.

\item{[Sw2]} J. \'Swi\c atkowski, {\it Hyperbolic Coxeter groups with
Sierpi\'nski carpet boundaries,} preprint 2015,  arXiv:1502.02044.

\item{[W]} D. Wise, From Riches to Raags: 3-manifolds, Right-Angled Artin Groups,
and Cubical Geometry, CBMS Regional Conference Series in Mathematics, Number 117,
AMS, 2012.

\bye